\documentclass[12pt,a4paper]{article}
\usepackage{amsmath}
\usepackage{amsfonts}
\usepackage[english,french]{babel}
\usepackage{fancyhdr}

\newcommand{\datum}{Version 99.07.22}
\date{June 1999}
\title{Equity Allocation and Portfolio Selection in Insurance: A simplified Portfolio Model}
\author{Erik Taflin\footnote{AXA,
23, Avenue Matignon, 75008 Paris, France; erik.taflin@u-bourgogne.fr, erik.taflin@axa.com}}

\newtheorem{theorem}{Theorem}[section]
\newtheorem{lemma}[theorem]{Lemma}
\newtheorem{proposition}[theorem]{Proposition}
\newtheorem{corollary}[theorem]{Corollary}

\numberwithin{equation}{section}

\begin{document}
\selectlanguage{english}

\maketitle
\thispagestyle{empty}
\begin{abstract}
\noindent
A quadratic discrete time probabilistic model, for optimal portfolio selection in
(re-)insurance is studied. For positive values of underwriting levels, the
expected value of the accumulated result is
optimized, under constraints on its variance and on annual ROE's. Existence of a
unique solution is proved and a Lagrangian formalism is
given. An effective method for solving the Euler-Lagrange equations is
developed. The approximate determination of the multipliers is discussed. This
basic model is an important building block for more complete models.  \\

\end{abstract}

\section{Introduction}
\label{Intro}
Optimal equity allocation and portfolio selection, for a reinsurance company with several portfolios
(or subsidiaries)
leads, in the presence of constraints on non-solvency probabilities,
market-shares and ROE's, to highly non-linear problems. This general situation
is studied in a forthcoming reference \cite{T2}. In certain situations
it is possible to replace the
non-solvency probabilities by stronger quadratic constraints (see \cite{T2}).
This gives a new simplified quadratic stochastic optimization problem, which solutions
(if any) respect the constraints and are approximate solutions of the original
problem. The purpose of the present paper is the resolution of the constraint quadratic
optimization problem in its most basic setting: \emph{Optimization of the expected
utility (i.e. the final accumulated result) of one single portfolio, under constraints on the
variance of the utility, on the annual ROE's and on the sign of the underwriting
levels, which should be positive.} This case covers the most basic
applications and is an important building block in more general situations \cite{T2}.
The portfolio in this paper is an extension of Markowitz portfolio \cite{M1}
to a multiperiod stochastic portfolio, as suggested by \cite{H-K} (c.f. also
\cite{D-J}). We here construct the portfolio such that
future results of contracts written at different times are distinguishable, which
easily allows to consider different maturity times.

We consider a portfolio
$\eta$ of $N$ different types of insurance contracts concluded at times $0, \ldots , \bar{T},$ where
$\bar{T} \geq 1.$ The amount of the contract
of type $i,$ where $1 \leq i \leq N,$ being concluded at time $t,$ where
$0 \leq t \leq  \bar{T},$ is denoted $\eta_{i}(t).$ In other words, $\eta_{i}(t)$
is the number of unit contracts (e.g. the unit is set to one FF)
of type $i.$ We suppose that the portfolio
$\eta$ does not generate any financial flows  at and after a certain time $\bar{T}+T.$ Let
$(\Delta U)(t+1,\eta)$ be the result of  the portfolio
$\eta$ for the period $[t,t+1[,$
let the utility of $\eta$, at time $t,$ $U(t,\eta)=\sum_{1 \leq s \leq t}(\Delta U)(s,\eta),$
($U(0,\eta)=0$), be the accumulated result for the period $[0,t[,$
and let the final utility of the portfolio
$\eta,$ $U(\infty,\eta)=U(\bar{T}+T,\eta)$ be the accumulated result
until no more financial flows are generated.
(For the precise definition of $U$ see formula (\ref{EqA.0.1}).)
For $1 \leq i \leq N$ and $0 \leq t \leq  \bar{T},$
$\eta_{i}(t)$ and $(\Delta U)(t+1,\eta)$ are random variables.
In the case under consideration the equity $K(t)=K(0)+U(t,\eta),$ where $K(0) \geq 0$ is the
initial equity at $t=0.$
Let $\mathcal{F}_{t}$
be the events which are possible up to time $t.$ The number of unit-contracts $\eta_{i}(t)$
to conclude at time $t,$ shall be known with certainty at time $t.$ This means
that $\eta_{i}(t)$ is $\mathcal{F}_{t}$--measurable. We impose
that the random variable $\eta_{i}(t)$ has finite variance, for $1 \leq i \leq N$
and $0 \leq t \leq  \bar{T}.$
We now introduce
the following constraints on the variable $\eta,$ ($E$ is the expectation operator):
\begin{itemize}
\item $C_{3}$) $E((\Delta U)(t+1,\eta))
				\geq c(t)E(K(t)),$ where
				$c(t) \in \mathbb{R}^{+}$ is given (constraint on profitability)
\item $C_{4}$) $E(( U(\infty, \eta)
				-E(U(\infty, \eta)))^{2})
				\leq \sigma^{2},$ where $\sigma^{2} > 0$ is given
				( acceptable level of the variance of the final utility)
\item $C_{6}$) $0 \leq \eta_{i}(t),$ where $1 \leq i \leq N$
				(only positive underwriting levels)
\end{itemize}
Let $\mathcal{C}_{0}$ be the set of portfolios $\eta,$
such that $\eta$ satisfies constraints ($C_{3}$), ($C_{4}$) and ($C_{6}$) and
such that $\eta_{i}(t)$
is $\mathcal{F}_{t}$--measurable and has finite second order moment, for $1 \leq i \leq N$
and $0 \leq t \leq  \bar{T}.$ We consider the problem of optimizing the expected
accumulated result, i.e. to find all
$\hat{\eta} \in \mathcal{C}_{0},$ such that
\begin{equation}
  E(U(\infty, \hat{\eta}))
   =\sup_{\eta \in \mathcal{C}_{0}} E(U(\infty, \eta)).
 \label{EqA.0}
\end{equation}
We establish (Corollary \ref{CorA.1}) under certain mild conditions, ($H_{1}$), ($H_{2}$) and ($H_{3}$) of \S 2,
on the result processes for the unit-contracts,
that the optimization problem (\ref{EqA.0}) has a unique solution, if $\mathcal{C}_{0}$
is non empty. Moreover the solution, is derived from  a Lagrangian formalism (\ref{EqA.13}),
and is given by formula (\ref{EqA.16}). Condition ($H_{1}$) says that the final
utility (sum of all results) of a unit contract, written at time $k,$
is independent of events occurred before $k.$ In practice, this is generally not
true, among other things, because of feed-forward phenomena in the prizing.
Condition ($H_{2}$) is equivalent to the statement that no non-trivial linear
combinations of final utilities, of contracts written at a certain time, is a
certain random-variable. This can also be coined, in more financial terms:
a underwriting portfolio $\eta(t),$ constituted at time $t,$ can not be risk-free.
Condition ($H_{3}$) says that the final
utility  of  unit contracts, written at different times are independent.
The conditions ($H_{1}$), ($H_{2}$) and ($H_{3}$), which
excludes many interesting situations, like cyclic markets, have only been chosen
for simplicity. They can largely be weakened, without altering the results of this paper.
An important point is that no particular distributions (statistical laws) are
required.

To determine, practically, the solution $\hat{\eta}$ given by formula (\ref{EqA.16}),
we shall give, in Appendix~\ref{Operator C}, an effective method to calculate the inverse  of the
linear operator representing the quadratic form in constraint ($C_{4}$). We
also establish that this operator has a finite spectrum. In
Appendix~\ref{Multipliers}, the determination of Lagrange multipliers in the solution
(\ref{EqA.16}) is studied and a simple approximation method is proposed. The proofs
of the mathematical
results, of \S~\ref{Section 2}, Appendix~\ref{Operator C} and Appendix~\ref{Multipliers}, are
given in Appendix~\ref{Proof}. This paper, which is self contained, is a formalized
version of the report \cite{T1}. \\

\noindent
{\bf Acknowledgement:} The author would like to thank Jean-Marie Nessi, CEO of
AXA-R\'e, and his collaborators for the many interesting discussions, which
were the starting point of this work.

\section{The model and main results}
\label{Section 2}

We define the portfolio in a probabilistic context, given by a probability
space  $(\Omega, P, \mathcal{F})$ and a filtration 
$\mathcal{A}=\{\mathcal{F}_{t}\}_{t \in \mathbb{N}},$ of sub $\sigma$-algebras of
the $\sigma$-algebra $\mathcal{F},$ i.e.
$\mathcal{F}_{0}=\{\Omega,\emptyset \}$ and
$\mathcal{F}_{s} \subset \mathcal{F}_{t} \subset \mathcal{F}$ for
$0 \leq s \leq t.$ The portfolio is composed by $N \geq 1$ types of insurance
contracts. The utility $u_{i}(t,t'),$ at $t' \in \mathbb{N}$ of the unit contract $i,$ 
$1 \leq i \leq N,$ concluded at time $t \in \{0,\ldots ,\bar{T} \}$ is by definition
the accumulated result in the time interval $[t,t'[$ for $t'>t,$
and $u_{i}(t,t')=0,$ for $0 \leq t' \leq t.$ We suppose that $u_{i}(t,t')$
is $\mathcal{F}_{t'}$-measurable and that $(u(t,t'))_{t' \geq 0}$ is an element
in the space\footnote{Let $1 \leq q < \infty.$ Then 
  $(X_{i})_{0 \leq i} \in \mathcal{E}^{q}(\mathbb{R}^{N})$
  if and only if  $X_{i}:\Omega \rightarrow \mathbb{R}^{N}$ is $\mathcal{F}$
  measurable and
  $\|X_{i}\|_{L^{q}(\Omega, \mathbb{R}^{N})}
  =(\int_{\Omega} |X_{i}(\omega)|^{q}_{\mathbb{R}^{N}}dP(\omega))^{1/q}
  < \infty$
  for $i \geq 0$, where $|\thickspace \thickspace|_{\mathbb{R}^{N}}$
  is the norm in $\mathbb{R}^{N}.$  Let $\mathcal{E}^{q}
  (\mathbb{R}^{N},\mathcal{A})$ the sub-space of
   $\mathcal{A}$ adapted processes in
  $\mathcal{E}^{q}(\mathbb{R}^{N}).$ We define $\mathcal{E}(\mathbb{R}^{N})
  =\cap_{q \geq 1}\mathcal{E}^{q}(\mathbb{R}^{N})$ and
  $\mathcal{E}(\mathbb{R}^{N},\mathcal{A})
  =\cap_{q \geq 1}\mathcal{E}^{q}(\mathbb{R}^{N},\mathcal{A}).$}
$\mathcal{E}(\mathbb{R}^{N}),$ of processes, with finite moments at all orders.
Since, for given $t \in \{0,\ldots ,\bar{T} \},$ the process $(u_{i}(t,t'))_{t' \geq 0}$
is $\mathcal{A}$-adapted, it follows that
$(u(t,t'))_{t' \geq 0} \in \mathcal{E}  (\mathbb{R}^{N},\mathcal{A}).$
The  final utility of the unit contract $i,$
$u_{i}^{\infty}(t)=u_{i}(t,s')$ $(=u_{i}(t,\infty)),$ when the contract
does not generate a  flow after the time $s',$ $s' \geq 0,$ is
$\mathcal{F}_{s'}$ measurable. By hypothesis (see \S 1) there exists here $s',$
such that $0 \leq s' \leq \bar{T}+T.$
We define the utility $U(t,\eta)$ of a portfolio $\eta$ at time $t \in \mathbb{N},$
where $\eta_{i}(s)$ is $\mathcal{F}_{s}$-measurable for $1 \leq i \leq N$ and
$s \in \{0,\ldots ,\bar{T} \},$ by\footnote{The scalar product in
  $\mathbb{R}^{N}$ is $x \cdot y =\sum_{1 \leq i \leq N} x_{i} y_{i}.$}
\begin{equation}
  U(t, \eta)
   =\sum_{k \leq t} \eta(k) \cdot u(k,t).
 \label{EqA.0.1}
\end{equation}
The result of $\eta$ for the period $[t,t+1[,$ introduced in \S 1, is
then given by $(\Delta U)(t+1,\eta)$ $=U(t+1, \eta)-U(t, \eta).$
The final utility $U(\infty,\eta)$ and the equity $K(t)$ are given as in \S 1.
Let $\mathcal{E}_{\bar{T}}^{q}  (\mathbb{R}^{N},\mathcal{A}),$ $1 \leq q$
(resp. $\mathcal{E}_{\bar{T}}  (\mathbb{R}^{N},\mathcal{A})$) be the subspace
of elements $\xi \in \mathcal{E}^{q}  (\mathbb{R}^{N},\mathcal{A})$
(resp. $\mathcal{E}  (\mathbb{R}^{N},\mathcal{A})$), with $\xi(t)=0$ for $t > \bar{T}.$
In optimization problem (\ref{EqA.0}) it is imposed that the portfolio has
finite variance, so it is an element of the Hilbert space
$\mathcal{H}= \mathcal{E}_{\bar{T}}^{2}  (\mathbb{R}^{N},\mathcal{A}),$ with
scalar product given by
$(\eta,\eta')_{\mathcal{H}}=\sum_{0 \leq t \leq \bar{T}} \int_{\Omega}
                 ((\eta (t))(\omega) \cdot (\eta' (t))(\omega))dP(\omega).$

The set $\mathcal{C}_{0}$ is well-defined, although the variance of the final
utility $U(\infty,\eta)$ is not finite in general, for the (unit-) utility processes
$(u(t,t'))_{t' \geq 0} \in \mathcal{E}  (\mathbb{R}^{N},\mathcal{A}).$ In fact the
quadratic form
\begin{equation}
 \eta \mapsto \mathfrak{a}(\eta)= E(( U(\infty,\eta))^{2}),
 \label{EqA.-1N}
 \end{equation}
in $\mathcal{H},$ has a maximal domain $\mathcal{D}(\mathfrak{a}),$
since for each $\eta \in \mathcal{H},$
the  stochastic process $(U(t, \eta))_{t \geq 0}$ is an element of the
space $\mathcal{E}^{p}(\mathbb{R}, \mathcal{A}),$ for $1 \leq p < 2$ (which
follows directly from Schwarz inequality).
The solution of optimization problem (\ref{EqA.0}) is largely based on the study of the
quadratic form
\begin{equation}
 \eta \mapsto \mathfrak{b}(\eta)= E(( U(\infty,\eta) -E(U(\infty,\eta)))^{2}),
 \label{EqA.1N}
 \end{equation}
in $\mathcal{H},$ with (maximal) domain
$\mathcal{D}(\mathfrak{b})=\mathcal{D}(\mathfrak{a}).$

We shall introduce another optimization problem, having only
piece-vice linear constraints and which we will prove, to be equivalent to
(\ref{EqA.0}), in a precise way.
Since $K(0) \geq 0,$ it follows from the formula (see \S\ref{Intro}) for the equity
$K(t)$ and constraint  $C_{3}$
that $E(U(\infty, \hat{\eta})) \geq 0,$ if the solution $\hat{\eta}$ exists.
Moreover, if $a \geq 1$ and constraint $C_{3}$ (resp. $C_{6}$) is satisfied,
then it is also satisfied with $\eta$ replaced by $a \eta.$ Therefore, if the
solution $\hat{\eta}$ exists, then
$E(( U(\infty, \hat{\eta}) -E(U(\infty, \hat{\eta})))^{2}) = \sigma^{2}.$
With this observation in mind, we introduce the optimization problem, which is
to find all $\hat{\eta} \in \mathcal{C}_{1}$ such that
\begin{equation}
\mathfrak{b}(\hat{\eta})
   =\inf_{\eta \in \mathcal{C}_{1}} \mathfrak{b}(\eta),
 \label{EqA.2N}
\end{equation}
where $\mathcal{C}_{1}$ is the set of all $\eta \in \mathcal{H}$ such that
the following constraints are satisfied:
\begin{itemize}
\item $C_{3}'$) $E((\Delta U)(t+1,\eta))
				\geq c(t)E(K(t)),$ where
				$c(t) \in \mathbb{R}^{+}$ is given (constraint on profitability)
\item $C_{4}'$) $E( U(\infty, \eta)) \geq e,$ where $e \geq 0$ is given
				(acceptable level of the expected final utility)
\item $C_{6}'$) $\eta_{i}(t) \geq 0,$ where $1 \leq i \leq N$
				(only positive underwriting levels).
\end{itemize}
$\mathcal{C}_{1}$ is a closed convex subset of $\mathcal{H}.$ In fact, the
constraint functions in $C_{3}'$ and $C_{4}'$ are strongly continuous, so it is
the intersection of  the three  closed convex cones defined by $C_{3}',$
$C_{4}'$ and $C_{6}'.$
As we will see, the advantage of this formulation is that it is easy to prove an
existence result.

We also introduce the auxiliary optimization problem, which consists of finding all
$\hat{\eta} \in \mathcal{C}_{2}$ such that
\begin{equation}
\mathfrak{a}(\hat{\eta})
   =\inf_{\eta \in \mathcal{C}_{2}} \mathfrak{a}(\eta),
 \label{EqA.3N}
\end{equation}
where $\mathcal{C}_{2}$ is the set of all $\eta \in \mathcal{H}$ such that
the following constraints are satisfied:
\begin{itemize}
\item $C_{3}''$) $E((\Delta U)(t+1,\eta))
				\geq c(t)E(K(t)),$ where
				$c(t) \in \mathbb{R}^{+}$ is given (constraint on profitability)
\item $C_{4}''$) $E( U(\infty, \eta)) = e,$ where $e \geq 0$ is given
				(acceptable level of the expected final utility)
\item $C_{6}''$) $\eta_{i}(t) \geq 0,$ where $1 \leq i \leq N$
				(only positive underwriting levels).
\end{itemize}
$\mathcal{C}_{2}$ is a closed convex subset of $\mathcal{H}.$

To solve the optimization problems (\ref{EqA.2N}) and (\ref{EqA.3N}), we first
establish properties of the positive quadratic forms $\mathfrak{b}$ and
$\mathfrak{a}.$ Let $\mathfrak{b}_{0}$ and $\mathfrak{a}_{0}$ be the restriction
of $\mathfrak{b}$ and $\mathfrak{a},$ respectively, to the domain
$\mathcal{D}(\mathfrak{b}_{0})=\mathcal{D}(\mathfrak{a}_{0})=
                  \mathcal{E}_{\bar{T}}(\mathbb{R}^{N},\mathcal{A}),$
where $\mathcal{E}_{\bar{T}}(\mathbb{R}^{N},\mathcal{A})
           =\mathcal{H} \cap \mathcal{E}(\mathbb{R}^{N},\mathcal{A}).$
The corresponding bilinear form is denoted by the same symbol 
$\mathfrak{b},$ $\mathfrak{a},$ $\mathfrak{b}_{0}$ and $\mathfrak{a}_{0},$
respectively,
i.e.
\begin{equation}
 \mathfrak{b}(\xi,\eta)=
   E(( U(\infty,\xi) -E(U(\infty,\xi)))( U(\infty,\eta) -E(U(\infty,\eta)))),
 \label{EqA.4N}
 \end{equation}
where $\xi, \eta \in \mathcal{D}(\mathfrak{b}),$ etc.
\begin{lemma}\label{LemA.2}
 $\mathfrak{b}_{0}(\eta),$ with domain
  $\mathcal{D}( \mathfrak{b}_{0})=\mathcal{E}_{\bar{T}}(\mathbb{R}^{N},\mathcal{A})$
(resp. $\mathfrak{a}_{0},$ with domain
  $\mathcal{D}( \mathfrak{a}_{0})=\mathcal{E}_{\bar{T}}(\mathbb{R}^{N},\mathcal{A})$)
is a densely defined closeable symmetric quadratic form in
 $\mathcal{H}.$
\end{lemma}
We make certain (technical) hypotheses on the claim processes:
\begin{itemize}
\item $H_{1}$) $u^{\infty}(k)$ is independent of $\mathcal{F}_{k}$ for
				$k \in \mathbb{N}$
\item $H_{2}$) for $k \in \mathbb{N}$ the  $N  \times N$ (positive) matrix
				$c(k)$ with elements
				$c_{ij}(k)=E((u_{i}^{\infty}(k)-E(u_{i}^{\infty}(k))
				(u_{j}^{\infty}(k)-E(u_{j}^{\infty}(k)))$ is strictly positive
\item $H_{3}$) $u_{i}^{\infty}(k)$ and $u_{j}^{\infty}(l)$ are independent for
				$k \neq l$
\end{itemize}
The first hypothesis implies boundedness of the quadratic forms.
\begin{lemma}\label{LemA.3}
If the hypothesis ($H_{1}$) is satisfied, then
$\mathcal{D}(\mathfrak{b})=\mathcal{D}(\mathfrak{a})=\mathcal{H}$
and  the quadratic forms $\mathfrak{b}$
and $\mathfrak{a}$ are bounded, i.e. there exists
$C^{2} > 0$ such that, 
\begin{equation}
  \mathfrak{b}(\eta) \leq
 C^{2} \| \eta \|^{2}_{\mathcal{H}} \quad \text{and} \quad
  \mathfrak{a}(\eta) \leq
 C^{2} \| \eta \|^{2}_{\mathcal{H}},
 \label{EqA.4}
 \end{equation}
for
$\eta \in \mathcal{H}.$
\end{lemma}
This lemma and the next crucial result show that (the square root of) each one of
the quadratic forms
$\mathfrak{b}$ and $\mathfrak{a}$ is equivalent to the norm in $\mathcal{H}.$
\begin{theorem}\label{ThA.4}
If the hypotheses ($H_{1}$), ($H_{2}$) and ($H_{3}$) are satisfied, then the
quadratic forms $\mathfrak{b}$ and $\mathfrak{a}$
are bounded from below by  strictly positive numbers $c_{1}^{2}$ respectively $c_{2}^{2},$
where $c_{1}^{2} \geq c_{2}^{2}>0,$  i.e.
\begin{equation}
  \mathfrak{b}(\eta) \geq c_{2}^{2} \| \eta \|^{2}_{\mathcal{H}} \quad \text{and} \quad
 \mathfrak{a}(\eta) \geq c_{1}^{2} \| \eta \|^{2}_{\mathcal{H}},
 \label{EqA.4.1}
 \end{equation}
for $\eta \in \mathcal{H}.$
\end{theorem}
Let $B$ (resp. $A$)  be the operator in
$\mathcal{H},$ associated with $\mathfrak{b}$ (resp. $\mathfrak{a}$),
(by the representation theorem),
i.e.
\begin{equation}
  \mathfrak{b}(\xi,\eta)=(\xi,B \eta)_{\mathcal{H}}
 \quad (\text{resp.} \quad
  \mathfrak{a}(\xi,\eta)=(\xi,A \eta)_{\mathcal{H}}),
 \label{EqA.5}
 \end{equation}
for $\xi \in \mathcal{H}$ and $\eta \in \mathcal{H}.$
$B$ and $A$ are strictly positive, bounded, self-adjoint, onto ${\mathcal{H}}$ and
invertible. The
inverse operators $B^{-1} : {\mathcal{H}} \rightarrow {\mathcal{H}}$ and
 $A^{-1} : {\mathcal{H}} \rightarrow {\mathcal{H}}$ are bounded
self-adjoint operators.  There exist $c \in \mathbb{R},$
such that $0 < cI \leq B \leq A,$ where $I$ is the identity operator.
It follows from formula (\ref{EqA.5}), that an explicit expression of $A$ is
given by
\begin{equation}
  (A \eta)(k)=E( U(\infty,\eta) u^{\infty}(k) | \mathcal{F}_{k})
 \label{EqA.6.1}
 \end{equation}
and that an explicit expression of $B$ is
given by
\begin{equation}
  (B \eta)(k)=E(( U(\infty,\eta) -E(U(\infty,\eta)))u^{\infty}(k) | \mathcal{F}_{k}),
 \label{EqA.6.2}
 \end{equation}
for $\eta \in \mathcal{H},$ where $0 \leq k \leq \bar{T}.$

\begin{theorem}\label{ThA.5}
Let the hypotheses ($H_{1}$), ($H_{2}$) and ($H_{3}$) be satisfied.
If $\mathcal{C}_{1}$ (resp.  $\mathcal{C}_{2}$) is nonempty, then there exists
a unique solution $\hat{\eta} \in \mathcal{C}_{1}$ (resp.  $\mathcal{C}_{2}$) of
the optimization problem for $\mathfrak{b}$ (resp. $\mathfrak{a}$)
\begin{equation}
\mathfrak{b}(\hat{\eta})
   =\inf_{\eta \in \mathcal{C}_{1}} \mathfrak{b}(\eta) \quad 
     (\text{resp.} \quad   \mathfrak{a}(\hat{\eta})
   =\inf_{\eta \in \mathcal{C}_{2}} \mathfrak{a}(\eta)).
 \label{EqA.12}
\end{equation}
\end{theorem}
Lemma \ref{LemA.3} shows that (the square root of) the function in
constraint ($C_{4}$) is strongly continuous in $\mathcal{H}$ and convex. So
$\mathcal{C}_{0}$ is a closed convex subset of $\mathcal{H}.$
\begin{lemma}\label{LemA.4}
Let the hypotheses ($H_{1}$), ($H_{2}$) and ($H_{3}$) be satisfied.
Solutions of optimization problems (\ref{EqA.0}) and (\ref{EqA.2N}) then have
the following properties: (i) If $\hat{\eta}$ is a solution of equation
(\ref{EqA.0}), then $\hat{\eta}$ is also a solution of equation (\ref{EqA.2N}),
with $E(U(\infty, \hat{\eta}))=e,$
(ii) If $\hat{\eta}$ is a solution of equation (\ref{EqA.2N}), then
$\hat{\eta}$ is also a solution of equation (\ref{EqA.0}), with
$\sigma^{2}=E(( U(\infty, \hat{\eta}) -E(U(\infty, \hat{\eta})))^{2}).$
\end{lemma}
This solves the original problem of \S\ref{Intro}:
\begin{corollary}\label{CorA.1}
Let hypotheses ($H_{1}$), ($H_{2}$) and ($H_{3}$) be satisfied. If
$\mathcal{C}_{0}$ is non-empty, then optimization problem (\ref{EqA.0}) has
unique solution $\hat{\eta} \in \mathcal{C}_{0}.$
\end{corollary}
To construct the solution of the optimization problem (\ref{EqA.12}), in the case of
$\mathfrak{b}$ (resp. $\mathfrak{a}$), we shall consider a Lagrangian
$h_{\lambda, \mu, \nu}^{\mathfrak{b}}$ (resp. $h_{\lambda, \mu, \nu}^{\mathfrak{a}}$),
with multipliers. Let $\lambda_{0}, \lambda_{1}, \ldots \lambda_{\bar{T}+T-1},$
and $\mu$ be real numbers and let $\nu \in \mathcal{H}.$
These are the multipliers.
A Lagrangian is defined by
\begin{equation}
  \begin{split}
  h_{\lambda, \mu, \nu}^{\mathfrak{c}} &( \eta) \\
              &=\frac{1}{2} \mathfrak{c}(\eta)
    -\sum_{0 \leq t \leq \bar{T}+T-1} \lambda_{t}(E((\Delta U)(t+1,\eta))- c(t)E(K(t))) \\
    &-\mu (E( U(\infty, \eta))-e) -(\nu,\eta)_{\mathcal{H}},
 \end{split} \label{EqA.13}
\end{equation}
where $\eta \in \mathcal{H},$ with $\mathfrak{c}=\mathfrak{b}$ (resp. $\mathfrak{a}$).  We now have to
find the critical points in $\mathcal{H},$ 
for fixed $\lambda,$ $\mu$ and $\nu$ and determine the multipliers
such that the critical point in fact is an element of $\mathcal{C}_{1}$ (resp. $\mathcal{C}_{2}$).
The multipliers shall satisfy
\begin{equation}
    \mu \geq 0, \; \mu (E( U(\infty, \hat{\eta}))-e)=0  \quad (\text{resp.} \; \mu  \in \mathbb{R}),
 \label{EqA.13.1}
\end{equation}
\begin{equation}
    \lambda_{t} \geq 0,
     \; \lambda_{t}(E((\Delta U)(t+1,\hat{\eta}))- c(t)E(K(t,\hat{\eta} )))=0,
 \label{EqA.14}
\end{equation}
for $0 \leq t \leq \bar{T}+T-1$ and
\begin{equation}
   \nu_{i}(k) \geq 0, \; (\nu_{i}(k))(\omega)=0 \quad \text{a.e. for}
           \quad \omega \in supp\, \hat{\eta}_{i}(k),
 \label{EqA.15}
\end{equation}
$0 \leq k \leq \bar{T}$ and $1 \leq i \leq N.$
Let $\mathcal{C}_{1}$ (resp. $\mathcal{C}_{2}$) be non-empty.
\begin{theorem}\label{ThA.7}
There are multipliers, satisfying (\ref{EqA.13.1}), (\ref{EqA.14})
and (\ref{EqA.15}), such that the solution $\hat{\eta},$ of the optimization
problem (\ref{EqA.12}) is given by the unique solution $\hat{\eta}$ of the
equation $(D h_{\lambda, \mu, \nu}^{\mathfrak{b}})(\eta )=0$
(resp. $(D h_{\lambda, \mu, \nu}^{\mathfrak{a}})(\eta )=0$).
Moreover,
\begin{equation}
   \hat{\eta}=C^{-1}(\mu \, m +\sum_{0 \leq t \leq \bar{T}+T-1} \lambda_{t} l_{t}
   +\nu),
 \label{EqA.16}
\end{equation}
where $C=B$ (resp. $C=A$) and where in this formula $m \in \mathcal{H}$  is
given by the linear functional $(m,\eta)_{\mathcal{H}}=E( U(\infty, \eta))$
and $l_{t} \in \mathcal{H}$ is given by the linear functional
$(l_{t},\eta)_{\mathcal{H}}=E((\Delta U)(t+1,\eta))- c(t)E(U(t,\eta)),$ for
$0 \leq t \leq \bar{T}+T-1.$
\end{theorem}
 Uniqueness
of the multipliers is assured if, for example, $m, l_{0}, \ldots ,l_{\bar{T}+T-1}$
are linearly independent and $c(t)>0,$ for $0 \leq t \leq \bar{T}+T-1.$
An explicit expression of $m$ is given by $m(k)=E( u^{\infty}(k) | \mathcal{F}_{k}).$
Hypothesis ($H_{1}$) then give
\begin{equation}
  m(k)=E( u^{\infty}(k)),
 \label{EqA.16.1}
\end{equation}
where $0 \leq k \leq \bar{T}.$
An explicit expression of $l_{t},$ $0 \leq t \leq \bar{T}+T-1,$ is given by
\begin{equation}
   l_{t}(k)=E( u(k,t+1) -(1+c(t))u(k,t) | \mathcal{F}_{k}),
 \label{EqA.16.2}
\end{equation}
if $0 \leq k \leq t$ and by $l_{t}(k)=0,$ if  $t< k,$
where in both cases $0 \leq k \leq \bar{T}.$ We note that, if it is supposed in
formula (\ref{EqA.16.2}), that $u(k,t)$ and $u(k,t+1)$ are independent of
$\mathcal{F}_{k},$ then
$l_{t}(k)=E( u(k,t+1) -(1+c(t))u(k,t)).$ So when this hypothesis, in addition to
hypothesis ($H_{1}$), is satisfied, then $l_{t}(k)$ is just a vector in
$\mathbb{R}^{N}.$

To determine, practically, the solution $\hat{\eta}$ given by formula
(\ref{EqA.16}), the operator $C^{-1}$ is needed. We shall give, in
Appendix~\ref{Operator C}, an effective method to calculate the inverse operator
$C^{-1}$ of $C,$ reminding that $C=A$ (resp. $C=B$) is defined by formula
(\ref{EqA.6.1}) (resp. (\ref{EqA.6.2})) in the case of $\mathfrak{a}$
(resp. $\mathfrak{b}$).

\appendix
\section{Appendix: The operator $C$}
\label{Operator C}
In this appendix we study the spectral properties of the operator $C,$
defined by  (\ref{EqA.6.1}) (resp. (\ref{EqA.6.2})) when $C=A$ (resp. $C=B$).
We also  give an algorithm, which determines $C^{-1}.$

The Hilbert space $\mathcal{H}$ has a canonical decomposition into a direct sum
$\mathcal{H}=\oplus_{0 \leq k \leq \bar{T}} \mathcal{H}_{k},$ where $\mathcal{H}_{k}$ is
the quotient space of $\mathcal{H}$ and the subspace of elements $\eta \in \mathcal{H},$
having $\eta (k) =0.$
To obtain the corresponding decomposition of the operators $A$ and $B,$
we introduce the notation
\begin{equation}
  A(k,l) \eta(l)=E( u^{\infty}(k) (u^{\infty}(l) \cdot \eta (l)) | \mathcal{F}_{k})
 \label{EqA.17}
 \end{equation}
and
\begin{equation}
  B(k,l) \eta(l)=E(u^{\infty}(k)(u^{\infty}(l) \cdot \eta (l)
          -E(u^{\infty}(l) \cdot \eta (l))) | \mathcal{F}_{k}),
 \label{EqA.18}
 \end{equation}
where $\eta \in \mathcal{H},$  $0 \leq k \leq \bar{T}$
and $0 \leq l \leq \bar{T}.$ Formula (\ref{EqA.6.1}) (resp. (\ref{EqA.6.2}))
then gives that $(A \eta)(k)=\sum_{0 \leq l \leq \bar{T}}A(k,l) \eta(l)$
(resp. $(B \eta)(k)=\sum_{0 \leq l \leq \bar{T}}B(k,l)$ $ \eta(l)$).
The continuity of the operators $A$ and $B$ in $\mathcal{H}$ shows that
$A(k,l)$ and $B(k,l)$ are continuous operators from $\mathcal{H}_{l}$ to
$\mathcal{H}_{k}.$

We also introduce the $N \times N$ real symmetric strictly positive matrices
$M^{\mathfrak{a}}(k)$ and $M^{\mathfrak{b}}(k)$, for $0 \leq k \leq \bar{T}.$ Here
$M_{ij}^{\mathfrak{a}}(k)=E(u_{i}^{\infty}(k) u_{j}^{\infty}(k))$
and
$M^{\mathfrak{b}}(k)$ have the same elements as the matrix $c(k)$ in hypothesis ($H_{2}$),
i.e. $M_{ij}^{\mathfrak{b}}(k)=E((u_{i}^{\infty}(k)-E(u_{i}^{\infty}(k))
               (u_{j}^{\infty}(k)-E(u_{j}^{\infty}(k))),$
for $0 \leq k \leq N$ and $0 \leq l \leq N.$
According to hypothesis ($H_{1}$), it follows that
\begin{equation}
  A(k,k) \eta(k)=M^{\mathfrak{a}}(k) \eta(k),
 \label{EqA.19}
 \end{equation}
for $0 \leq k \leq \bar{T}.$
The spectrum
$\sigma (M^{\mathfrak{a}}(k))$ (resp. $\sigma (M^{\mathfrak{b}}(k))$) of the matrix
$M^{\mathfrak{a}}(k)$ (resp. $M^{\mathfrak{b}}(k)$) is a finite set of strictly positive
real numbers, according to hypothesis ($H_{2}$). Let $\lambda \in \mathbb{R}$ be in
the resolvent set of $M^{\mathfrak{a}}(k).$
The inverse $(A(k,k)-\lambda)^{-1}$ is then simply obtained by inverting
the matrix $M^{\mathfrak{a}}(k)-\lambda.$
\begin{proposition}\label{PropA.1}
Let $\lambda \in \mathbb{R}$ be in
the resolvent sets of $M^{\mathfrak{a}}(k)$ and $M^{\mathfrak{b}}(k).$
Then the inverse
$(B(k,k)-\lambda)^{-1}$ of the operator $B(k,k)-\lambda$
exists and is given by
\begin{equation}
  \begin{split}
  (B(k,k)-\lambda)^{-1} \eta(k)&=(M^{\mathfrak{a}}(k)-\lambda)^{-1}
  ( \eta(k)-E(\eta(k))) \\
   + &(M^{\mathfrak{b}}(k)-\lambda)^{-1} E(\eta(k)),
 \end{split}  \label{EqA.20}
 \end{equation}
for $0 \leq k \leq \bar{T}.$
\end{proposition}
Let $C(k,l)=A(k,l)$ (resp. $B(k,l)$), for $C=A$ (resp. $B$).
To solve the equation $C \eta =\xi,$
where $\xi \in \mathcal{H}$ is given, we shall study
the system of equations
\begin{equation}
  \sum_{ 0 \leq l \leq \bar{T}} C(k,l) \eta(l)= \xi(k), \quad 0 \leq k \leq \bar{T}.
 \label{EqA.21}
 \end{equation}
 The operator $C(k,l)$ is continuous from $\mathcal{H}_{l}$
to $\mathcal{H}_{k}.$
Although, the operators $C(k,l)$ are integral operators in general, the system
(\ref{EqA.21}) can be solved in a finite number of steps only involving linear
algebra in $\mathbb{R}^{N}$ and evaluation of expectation values. To study also
the spectral properties of the operator $C$ in $\mathcal{H},$ we give this
algorithm for a more general equation in $\mathcal{H},$ namely
$C \eta -\lambda \eta =\xi,$ i.e.
\begin{equation}
  \sum_{ 0 \leq l \leq \bar{T}} (C(k,l) -\delta_{kl} \lambda) \eta(l)= \xi(k),
       \quad \xi(k) \in \mathcal{H}_{k}  \quad 0 \leq k \leq \bar{T},
 \label{EqA.22}
 \end{equation}
when $\lambda$ is in a certain subset, to be determined, of the resolventset
$\mathbb{R}-\sigma(C),$ of the operator $C.$

The first part of the algorithm consists of transforming the system (\ref{EqA.22})
into an equivalent lower triangular system.
For a given real number $\lambda$
outside a certain finite set (of poles),
we shall define operators $C^{n}(k,l)$ and
elements $\xi^{n} \in \mathcal{H},$ where $k,l,n \in \{0, \ldots ,\bar{T} \}.$
The definition will be specified by an finite iteration, beginning with $n=\bar{T}$
and ending with $n=0.$
Let us introduce the $N \times N$ real symmetric matrices $N_{n}^{\mathfrak{a}}(k)$ and
$N_{n}^{\mathfrak{b}}(k)$ and $D_{\bar{T}}(k))$ by
\begin{equation}
  (N_{n}^{\mathfrak{a}}(k))_{ij}
=E(E(u^{\infty}_{i}(k)| \mathcal{F}_{n}) E(u^{\infty}_{j}(k)| \mathcal{F}_{n}))
 \label{EqA.22.1}
 \end{equation}
and
\begin{equation}
  (N_{n}^{\mathfrak{b}}(k))_{ij} =E(
(E(u^{\infty}_{i}(k)| \mathcal{F}_{n}) - E(u^{\infty}_{i}(k)))
      (E(u^{\infty}_{j}(k)| \mathcal{F}_{n}) - E(u^{\infty}_{j}(k)))  ),
 \label{EqA.22.2}
\end{equation}
for $\quad 0 \leq i,j \leq N,$  $\quad 0 \leq k \leq \bar{T}$ and
$\quad 0 \leq n \leq \bar{T}.$
Symmetric matrices $D_{\bar{T}}(k),$ $0 \leq k \leq \bar{T}$
and real numbers $f_{\bar{T}},$ $g_{\bar{T}},$
and $d_{\bar{T}}$
are given by
\begin{equation}
  D_{\bar{T}}(k) = M^{\mathfrak{a}}(k)-\lambda, \quad f_{\bar{T}}=1, \quad g_{\bar{T}}=0,
 \label{EqA.22.3}
\end{equation}
and
\begin{equation}
  d_{\bar{T}} =
   E(u^{\infty}(\bar{T})) \cdot ((D_{\bar{T}}(\bar{T}))^{-1}E(u^{\infty}(\bar{T}))),
 \label{EqA.23}
\end{equation}
for $\lambda \in \mathbb{R}-\sigma(M^{\mathfrak{a}}(k)).$
For given integer $n,$ where $0 < n \leq \bar{T},$ given real numbers $f_{n},$ $g_{n}$
and given  $N \times N$ real symmetric matrices $D_{n}(k),$  where
$0 \leq k \leq n$ and where $D_{n}(n)$ is invertible, we define real numbers
$f_{n-1},$ $g_{n-1},$ $d_{n}$ and $N \times N$ real symmetric matrices $D_{n-1}(k):$
\begin{equation}
  d_{n} =
   E(u^{\infty}(n)) \cdot ((D_{n}(n))^{-1}E(u^{\infty}(n))),
 \label{EqA.23.1}
\end{equation}
\begin{equation}
   f_{n-1}=f_{n}(1-d_{n}f_{n}), \quad g_{n-1}=g_{n}+d_{n}f_{n}^{2},
 \label{EqA.23.2}
\end{equation}
and
\begin{equation}
 D_{n-1}(k) = D_{n}(k)-d_{n} f_{n}^{2} N_{n}^{\mathfrak{a}}(k), \quad 0 \leq k \leq n-1.
 \label{EqA.23.3}
\end{equation}
Let $d_{0}$ be given by formula (\ref{EqA.23.1}), with $n=0.$ The elements
$f_{n},$ $g_{n}$ and $d_{n}$ of the sequences $f_{0}, \ldots ,f_{\bar{T}},$
$g_{0}, \ldots ,g_{\bar{T}}$ and $d_{0}, \ldots ,d_{\bar{T}}$ are real valued
rational functions of $\lambda$ in $\mathbb{R}.$ Similarly, each matrix element
of the matrix $D_{n}(n),$ of the sequences $D_{0}(0), \ldots ,D_{\bar{T}}(\bar{T}),$
is a real valued rational function of $\lambda$ in $\mathbb{R}.$
We can now introduce the linear continuous operators
$C^{n}(k,l):\mathcal{H}_{l} \rightarrow \mathcal{H}_{k},$
for $k,l,n \in \{0, \ldots ,\bar{T} \}.$
We define $C^{\bar{T}}(k,l)$ by
\begin{equation}
 C^{\bar{T}}(k,l)=C(k,l), \quad  k \neq l
 \label{EqA.24.1}
\end{equation}
and
\begin{equation}
 C^{\bar{T}}(k,k)=C(k,k) -\lambda I,
 \label{EqA.24.2}
\end{equation}
(where $I$ is the identity operator on ${H}_{k}$),
for $k,l \in \{0, \ldots ,\bar{T} \}$ and $\lambda \in \mathbb{R}.$
For given integer $n,$ where $0 < n \leq \bar{T},$
and given  linear continuous operators
$C^{n}(k,l):\mathcal{H}_{l} \rightarrow \mathcal{H}_{k},$
where $k,l \in \{0, \ldots ,\bar{T} \}$ we define operators
$C^{n-1}(k,l):\mathcal{H}_{l} \rightarrow \mathcal{H}_{k}.$
If $\lambda \in \mathbb{R}$ is such that, $D_{n}(n)$ is invertible, then we define
\begin{equation}
 C^{n-1}(k,l)=C^{n}(k,l), \quad  n \leq k \leq \bar{T},
      \;  0 \leq l \leq \bar{T},
 \label{EqA.24.3}
\end{equation}
\begin{equation}
 C^{n-1}(k,l)=0, \quad 0 \leq  k < n, \; n \leq l \leq \bar{T},
 \label{EqA.24.4}
\end{equation}
\begin{equation}
 C^{n-1}(k,l)= f_{n-1}C(k,l),
   \quad  k \neq l, \; 0 \leq  k < n, \; 0 \leq  l < n,
 \label{EqA.24.5}
\end{equation}
\begin{equation}
 A^{n-1}(k,k) \eta(k)= D_{n-1}(k) \eta(k),
   \quad 0 \leq  k < n, \; \eta(k) \in \mathcal{H}_{k}
 \label{EqA.24.6}
\end{equation}
and
\begin{equation}
  \begin{split}
 B^{n-1}(k,k) \eta(k)&=D_{n-1}(k) \eta(k) \\
   &-(1-g_{n-1}) E(u^{\infty}(k)) (E(u^{\infty}(k)) \cdot E(\eta(k))), \\
   &\quad 0 \leq  k < n, \; \eta(k) \in \mathcal{H}_{k}.
 \end{split}  \label{EqA.24.7}
\end{equation}
If $C=A$ (resp. $B$), then $C^{n-1}(k,k)=A^{n-1}(k,k)$
(resp. $C^{n-1}(k,k)=B^{n-1}(k,k)$), for $0 \leq  k < n.$

According to (\ref{EqA.24.6}), the operator $A^{n}(k,k),$ in $\mathcal{H}_{k},$
has a bounded inverse, if $\lambda$ is such that, $D_{n}(k)$ is invertible:
\begin{equation}
 (A^{n}(k,k))^{-1} \eta(k)= (D_{n}(k))^{-1} \eta(k),
   \quad 0 \leq  k \leq n, \; \eta(k) \in \mathcal{H}_{k}.
 \label{EqA.25.1}
\end{equation}
As in
the case of $(B(k,k)-\lambda)^{-1},$ (see formula (\ref{EqA.20})), an explicit
expression for the inverse of $B^{n}(k,k)$ can be given.
\begin{proposition}\label{PropA.2}
Let $\lambda$ be
in the resolventset of the two operators
$M^{\mathfrak{c}}(k)
 -\sum_{n+1 \leq r \leq \bar{T}} d_{r} f_{r}^{2} N_{r}^{\mathfrak{c}}(k),$
for $\mathfrak{c}=\mathfrak{a}$ and $\mathfrak{c}=\mathfrak{b}.$ Then
\begin{equation}
  \begin{split}
  (B^{n}&(k,k))^{-1} \eta(k) \\
        &=(M^{\mathfrak{a}}(k)
     -\sum_{n+1 \leq r \leq \bar{T}} d_{r} f_{r}^{2} N_{r}^{\mathfrak{a}}(k) - \lambda)^{-1}
  ( \eta(k) -E(\eta(k))) \\
       &+(M^{\mathfrak{b}}(k)
        -\sum_{n+1 \leq r \leq  \bar{T}}  d_{r} f_{r}^{2} N_{r}^{\mathfrak{b}}(k) - \lambda)^{-1}
               E(\eta(k)),
 \end{split}  \label{EqA.25.2}
 \end{equation}
where  $0 \leq k \leq n \leq  \bar{T}.$
\end{proposition}
For reference we note that
\[ D_{n}(k)= M^{\mathfrak{a}}(k)
     -\sum_{n+1 \leq r \leq \bar{T}} d_{r} f_{r}^{2} N_{r}^{\mathfrak{a}}(k) - \lambda, \]
according to (\ref{EqP.9.2}), of Appendix~\ref{Proof}.

Next we define $\xi^{n-1} \in \mathcal{H}.$ For given $\xi \in \mathcal{H},$ let
$\xi^{\bar{T}}=\xi.$ For given integer $n,$ $0 < n \leq \bar{T},$ for given
linear continuous operators $C^{n}(k,l): \mathcal{H}_{l} \rightarrow \mathcal{H}_{k},$
$k,l \in \{0, \ldots ,\bar{T} \},$ where $C^{n}(n,n)$ has a bounded inverse, and
for given $\xi^{n},$ we define
\begin{equation}
  \xi^{n-1}(k)=\xi^{n}(k), \quad n \leq  k  \leq \bar{T}.
 \label{EqA.26}
\end{equation}
and
\begin{equation}
  \xi^{n-1}(k)=\xi^{n}(k) -C^{n}(k,n)(C^{n}(n,n))^{-1}\xi^{n}(n),
  \quad   0 \leq k < n.
 \label{EqA.27}
\end{equation}

To give spectral properties of the operator $C,$ we introduce notations for
unions of spectra of certain matrices. We recall that, if X is a linear operator
in a Hilbert space, then $\sigma (X)$ denotes here its spectrum. For commodity we
first introduce $\sigma_{\bar{T}+1}^{\mathfrak{a}}$ $=\sigma_{\bar{T}+1}^{\mathfrak{b}}$
$=\emptyset,$ the empty set. For $0  \leq n \leq  \bar{T},$ let
\begin{equation}
 \sigma_{n}^{\mathfrak{a}}=\sigma_{n+1}^{\mathfrak{a}} \cup
        \sigma(M^{\mathfrak{a}}(n)
 -\sum_{n+1 \leq r \leq \bar{T}} d_{r} f_{r}^{2} N_{r}^{\mathfrak{a}}(n))
 \label{EqA.27.1}
\end{equation}
and let
\begin{equation}
 \sigma_{n}^{\mathfrak{b}}=\sigma_{n}^{\mathfrak{a}} \cup
        \sigma(M^{\mathfrak{b}}(n)
 -\sum_{n+1 \leq r \leq \bar{T}} d_{r} f_{r}^{2} N_{r}^{\mathfrak{b}}(n)).
 \label{EqA.27.2}
\end{equation}
We note that $\sigma_{n}^{\mathfrak{c}} \subset \sigma_{n-1}^{\mathfrak{c}},$ for
$\mathfrak{c}=\mathfrak{a}$ and for $\mathfrak{c}=\mathfrak{b}.$
\begin{proposition}\label{PropA.3}
If $\lambda \notin \sigma_{n+1}^{\mathfrak{c}},$ then $C^{\bar{T}},\ldots ,C^{n}$
and $\xi^{\bar{T}},\ldots ,\xi^{n},$ given by (\ref{EqA.22.3})-(\ref{EqA.27})
are sequences of linear continuous operators in $\mathcal{H}$ and elements
of $\mathcal{H},$ respectively and
the system (\ref{EqA.22}) has the same solutions $\eta \in \mathcal{H}$ as the
system
\begin{equation}
  \sum_{ 0 \leq l \leq \bar{T}} C^{n}(k,l) \eta(l)= \xi^{n}(k),
                     \quad 0 \leq k \leq \bar{T}.
 \label{EqA.28}
\end{equation}
\end{proposition}
For $n=0,$ the system (\ref{EqA.28}) is lower triangular and reads
\begin{equation}
  \sum_{ 0 \leq l \leq k} C^{0}(k,l)  \eta(l)= \xi^{0}(k),
                     \quad 0 \leq k \leq \bar{T},
 \label{EqA.29}
\end{equation}
where we remind that $C^{0}$ and $\xi^{0}$ are dependent of
$\lambda \notin  \sigma_{1}^{\mathfrak{c}}.$

The second  part of the algorithm consists of solving the system (\ref{EqA.29})
by successive substitutions.
\begin{proposition}\label{PropA.4}
If $\lambda \notin  \sigma_{0}^{\mathfrak{a}}$
(resp. $\lambda \notin  \sigma_{0}^{\mathfrak{b}}$),
when $C=A$ (resp. $C=B$),
then $\eta(0), \ldots, \eta(\bar{T}),$ are successively given by
\begin{equation}
   \eta(k) =(C^{k}(k,k))^{-1}
   (\xi^{0}(k) -\sum_{ 0 \leq l \leq k-1} C^{0}(k,l) \eta(l)),
                     \quad 0 \leq k \leq \bar{T},
 \label{EqA.30}
\end{equation}
where we note that the sum is absent for $k=0.$
\end{proposition}
We remind that the operator $(C^{k}(k,k))^{-1},$ $0 \leq k \leq \bar{T},$
is explicitly given by (\ref{EqA.25.1}) (resp.  (\ref{EqA.25.2})), when $C=A$
(resp. $C=B$).

To sum up,
if $\lambda \notin  \sigma_{0}^{\mathfrak{a}}$
(resp. $\lambda \notin  \sigma_{0}^{\mathfrak{b}}$),
then $\eta,$ given by formula (\ref{EqA.30}), is the unique
solution, in $\mathcal{H},$ of equation (\ref{EqA.22}), in the case of $C=A$ (resp. $C=B$).
Consequently $\sigma(A) \subset \sigma_{0}^{\mathfrak{a}}$
and $\sigma(B) \subset \sigma_{0}^{\mathfrak{b}}.$
In particular the spectrum $\sigma(C)$ is a finite set of strictly positive
real numbers.

\section{The Lagrange multipliers}
\label{Multipliers}
We have already established that there exists multipliers $\mu,$
$\lambda_{0}, \ldots, \lambda_{\bar{T}+T-1} \in \mathbb{R}$ and $\nu \in \mathcal{H},$
such that the unique solution $\hat{\eta}$ of the optimization problem (\ref{EqA.12})
is given by formula (\ref{EqA.16}). In this appendix, we derive equations of which
the multipliers are solutions. These equations are taken as a starting point for
the derivation of algorithms permitting to calculate approximations of the
multipliers. Starting with the zeroth approximation as being the value of the multipliers
corresponding to the solution of the deterministic approximation of the optimization
problem (\ref{EqA.12}), we derive an explicit first approximation of the multipliers
in the stochastic case. More precise approximations is a subject of future studies.

A basic building block in our equations for the multipliers, is  the solution of
a certain finite dimensional optimization problem. (The dimension will mostly
be $N$ or $\bar{T}+T,$ i.e. the number of types of contracts or number of constraints
in $(C_{3}'$) and $(C_{4}'$) of \S~\ref{Section 2} respectively). For $z \in \mathbb{R}^{n},$
let $z \geq 0$ be defined by $z_{i} \geq 0$ for $1 \leq i \leq n.$ Let $m$ be a
strictly positive symmetric operator in $\mathbb{R}^{n}.$ For given $x \in \mathbb{R}^{n},$
let $F^{+}_{m}(x) \geq 0$ be the unique solution of the optimization problem
\begin{equation}
 \frac{1}{2} F^{+}_{m}(x) \cdot (mF^{+}_{m}(x)) -x \cdot F^{+}_{m}(x))
    =\inf_{y \geq 0} ( \frac{1}{2} y \cdot (my) -x \cdot y).
 \label{EqLM.1}
\end{equation}
The solution $F^{+}_{m}(x)$ of (\ref{EqLM.1}) is also the unique critical point of the
Lagrangian $L_{x}:$
\begin{equation}
 L_{x}(y)= \frac{1}{2} y \cdot (my) -x \cdot y -F^{-}_{m}(x) \cdot y,
 \label{EqLM.2}
\end{equation}
where $(F^{-}_{m}(x))_{i},$ $1 \leq i \leq n,$ are the multipliers corresponding
to the constraint $y \geq 0.$ We have
\begin{equation}
 F^{+}_{m}(x)= m^{-1} (x + F^{-}_{m}(x)),
 \label{EqLM.3}
\end{equation}
where
\begin{equation}
 F^{+}_{m}(x) \geq 0, \quad  F^{-}_{m}(x) \geq 0, \quad F^{+}_{m}(x) \cdot F^{-}_{m}(x)=0.
 \label{EqLM.4}
\end{equation}
We note that there exists, for given $m,$ a function $x \mapsto P_{x}$ of
$\mathbb{R}^{n}$ to linear projections on $\mathbb{R}^{n},$ (i.e. $P_{x}=P_{x}^{2}$),
taking a finite number ($2^{n}$) of values, such that $F^{+}_{m}(x)= m^{-1}P_{x}x.$

We shall here consider the solution $\hat{\eta }$ (given by (\ref{EqA.16})) only
in the case of $C=B.$ The case $C=A,$ is similar but simpler.
Let $\lambda_{\bar{T}+T}=\mu$ and let $l_{\bar{T}+T}=m.$ We define $\Theta^{\epsilon}(k)$
by
\begin{equation}
 \Theta^{\epsilon}(k)=F^{\epsilon}_{M^{\frak{a}}(k)}(
   \sum_{0 \leq t \leq \bar{T}+T} \lambda_{t} l_{t}(k)
       +(M^{\frak{a}}(k)-M^{\frak{b}}(k))E(\hat{\eta }(k)))
        -\sum_{l \neq k} B(k,l)\hat{\eta}(l)), 
 \label{EqLM.5}
\end{equation}
for $0 \leq k \leq \bar{T}$ and $\epsilon=\pm.$ It follows that
\begin{equation}
\hat{\eta}(k)=\Theta^{+}(k) \quad \text{and} \quad \nu(k)=\Theta^{-}(k),
 \label{EqLM.6}
\end{equation}
for $0 \leq k \leq \bar{T}.$

For simplicity we suppose
in the sequel of this paragraph that the subset $\{ E(l_{0}), \ldots ,E(l_{\bar{T} +T}) \}$
of $\mathcal{H}$ is linearly independent. This ensures that the multipliers are
unique in the deterministic case as well as in the stochastic case.
Expressions (\ref{EqLM.5}) and (\ref{EqLM.6}) are satisfied by the multipliers and
the solution $\hat{\eta}$ given by formula (\ref{EqA.16}).
To obtain another such equation
let $L$ be  the strictly positive symmetric operator in $\mathbb{R}^{\bar{T}+T}$ defined by
$(L^{-1})_{ts}=(l_{t},B^{-1}l_{s})$ and let $r(\nu) \in \mathbb{R}^{\bar{T}+T}$
be given by $(r(\nu))_{t}=(l_{t},B^{-1}\nu).$ Then the multiplier $\lambda$
satisfies:
\begin{equation}
\lambda=F^{-}_{L}(L(r(\nu)-e)).
 \label{EqLM.7}
\end{equation}
Equations (\ref{EqLM.5}), (\ref{EqLM.6}) and (\ref{EqLM.7}) form a closed system
for the solution $\hat{\eta}$ and the multipliers $\lambda$ and $\nu,$
(reminding that $\lambda_{\bar{T}+T}=\mu$). Alternatively, expression (\ref{EqA.16})
of $\hat{\eta},$ formula (\ref{EqLM.5}) with $\epsilon = -,$ the expression of $\nu$
in (\ref{EqLM.6})  and expression (\ref{EqLM.7}) of $\lambda$, also form a
closed system of equations, for $\hat{\eta},$ $\lambda$ and $\nu.$

We shall next give a crude approximation of the multipliers $\lambda$ and $\nu.$
To construct the zeroth approximation, we consider the deterministic approximation
of the optimization problem (\ref{EqA.12}). Let $\mathcal{C}_{D}$ be the subset
of elements $\eta \in \mathcal{C}_{1},$ such that $\eta(k)$ is $\mathcal{F}_{0}$
measurable for $0 \leq k \leq \bar{T}.$ The deterministic approximation, of the optimization
problem, is then to find all $\hat{\eta}_{D} \in \mathcal{C}_{D}$ such that
\begin{equation}
\frak{b}(\hat{\eta}_{D})
   =\inf_{\eta \in \mathcal{C}_{D}} \frak{b}(\eta).
 \label{EqLM.8}
\end{equation}
If $\mathcal{C}_{D}$ is non-empty, then there exists a unique solution $\hat{\eta}_{D}$
given by
\begin{equation}
   \hat{\eta}_{D}(k)
=(M^{\frak{a}}(k))^{-1}( \sum_{0 \leq t \leq \bar{T}+T} \lambda^{D}_{t} l^{D}_{t}(k) +\nu^{D}(k)),
 \label{EqLM.9}
\end{equation}
where $(l^{D}_{t}(k))_{i}=E((l_{t}(k))_{i})$ and where the multipliers
$\lambda^{D}_{t},$ $0 \leq t \leq \bar{T}+T$ and $\nu^{D}_{i}(k),$ $0 \leq k \leq \bar{T},$
$1 \leq i \leq N$ are real numbers such that
\begin{equation}
\lambda^{D}_{t} \geq 0, \quad (l^{D}_{t}, \hat{\eta}_{D})_{H} -e_{t} \geq 0,
   \quad \lambda^{D}_{t}((l^{D}_{t}, \hat{\eta}_{D})_{H} -e_{t})=0,
 \label{EqLM.10}
\end{equation}
for $0 \leq t \leq \bar{T}+T$ and
\begin{equation}
\nu^{D} \geq 0, \quad  \hat{\eta}_{D} \geq 0, \quad \nu^{D}_{i}(k)(\hat{\eta}_{D}(k))_{i}=0,
 \label{EqLM.11}
\end{equation}
for $0 \leq k \leq \bar{T}$ and $1 \leq i \leq N.$ The value of the multipliers
$\lambda^{D}$ and $\nu^{D}$ can be determined using the functions $F_{m}.$

Next we determine a first approximation $\lambda^{(1)}$ and $\nu^{(1)}$ of the
multipliers. Let
\begin{equation}
\lambda^{(1)}=F^{-}_{L}(L(r(\nu^{D})-e)),
 \label{EqLM.12}
\end{equation}
let
\begin{equation}
\bar{\eta}^{(1)}
=B^{-1}( \sum_{0 \leq t \leq \bar{T}+T} \lambda^{(1)}_{t} l_{t} +\nu^{D}),
 \label{EqLM.13}
\end{equation}
and let
\begin{equation}
\begin{split}
 \nu^{(1)}(k)&=F^{-}_{M^{\frak{a}}(k)}(
   \sum_{0 \leq t \leq \bar{T}+T} \lambda^{(1)}_{t} l_{t}(k) \\
       &+(M^{\frak{a}}(k)-M^{\frak{b}}(k))E(\bar{\eta }^{(1)}(k)))
        -\sum_{l \neq k} B(k,l)\bar{\eta }^{(1)}(l)), 
 \end{split} \label{EqLM.14}
\end{equation}
for $0 \leq k \leq \bar{T}.$ The corresponding first approximation $\hat{\eta}^{(1)}$
of the solution $\hat{\eta}$ is then defined by
\begin{equation}
\begin{split}
 \hat{\eta}^{(1)}(k)&=F^{+}_{M^{\frak{a}}(k)}(
   \sum_{0 \leq t \leq \bar{T}+T} \lambda^{(1)}_{t} l_{t}(k) \\
       &+(M^{\frak{a}}(k)-M^{\frak{b}}(k))E(\bar{\eta }^{(1)}(k)))
        -\sum_{l \neq k} B(k,l)\bar{\eta }^{(1)}(l)), 
 \end{split} \label{EqLM.15}
\end{equation}
for $0 \leq k \leq \bar{T}.$ Formula (\ref{EqLM.15}) ensures that $\hat{\eta}^{(1)} \geq 0.$
The approximation method defined by formulas (\ref{EqLM.13}), (\ref{EqLM.14})
and (\ref{EqLM.15}) can in an obvious way be generalized to higher order approximations.
However the convergence of the method must be established.

\section{Proofs}
\label{Proof}
We first give the proofs of the results of \S \ref{Section 2} and then those of
Appendix \ref{Operator C}.
\subsection{Proofs of results in \S \ref{Section 2}}
\label{Proof Section 2}
\noindent
\textbf{Proof of Lemma \ref{LemA.2}}
Schwarz inequality shows that $\eta \mapsto E(U(\infty,\eta))$ is a bounded
linear map from $\mathcal{H}$ to $L^{2}=L^{2}(\Omega, \mathbb{R}),$ so it is
enough to prove the statement in the case of $\mathfrak{a}.$ Since
$\mathfrak{a}$ is maximal, we only have to prove that the operator
$\mathcal{D}(\mathfrak{a}) \ni \eta \mapsto U(\infty,\eta) \in L^{2}$ in
$\mathcal{H},$ is a closeable, where
$L^{2}=L^{2}(\Omega, \mathbb{R}).$ Let $\{s_{n}\}_{n \geq 1}$ be a sequence
in $\mathcal{D}( \mathfrak{a}),$ such that
$\|s_{n}\|_{\mathcal{H}} \rightarrow 0$ and
$\|U(\infty,s_{n}) - v\|_{L^{2}} \rightarrow 0,$ when
$n \rightarrow \infty$ for some $v \in L^{2}.$ For $a \in \mathbb{R},$ let
$g_{a}$ be the characteristic function of the set
$\{\omega \in \Omega | (\sum_{0 \leq k \leq \bar{T}} |(u^{\infty}(k))(\omega)|^{2})^{1/2}
\leq a \}.$ Then
$E((g_{a}U(\infty,s_{n}))^{2})$ $\leq E(a^{2} \sum_{0 \leq k \leq \bar{T}}
|s_{n}(k)|^{2})$ $= a^{2} \|s_{n}\|^{2}_{\mathcal{H}}.$ Hence
$g_{a}U(\infty,s_{n})$ converges to $0$ in $ L^{2},$ when $n \rightarrow \infty$
for every $a \in \mathbb{R}.$ This gives that $v=0,$ which proves the statement. \\

\noindent
\textbf{Proof of Lemma \ref{LemA.3}}
 We prove that ($H_{1}$) implies that $\mathfrak{a}_{0}$ is bounded.
 Let $X_{k}=\eta (k) \cdot u^{\infty}(k),$ for $0 \leq k \leq \bar{T}.$ Then
 $U(\infty,\eta)=\sum_{0 \leq k \leq \bar{T}}X_{k},$ so
 $(E((U(\infty,\eta))^{2}))^{1/2} \leq
 \sum_{0 \leq k  \leq \bar{T}}(E((X_{k})^{2}))^{1/2}.$
 Since $\eta(k)$ is $\mathcal{F}_{k}$ measurable, since
 $E( |u^{\infty}_{i}(k)|^{2}  | \mathcal{F}_{k})$
 $=E( |u^{\infty}_{i}(k)|^{2})$
 according to hypothesis ($H_{1}$) and since
 $(X_{k})^{2}) \leq |u^{\infty}(k)|^{2} |\eta(k)|^{2}$
 it follows that
 $E((X_{k})^{2})=E(E((X_{k})^{2}| \mathcal{F}_{k} )) \leq
 E( |u^{\infty}(k)|^{2}) E(|\eta(k)|^{2}).$
 This shows that $\mathfrak{a}_{0}(\eta)
 =E((U(\infty,\eta))^{2}) \leq C^{2} \| \eta \|^{2}_{\mathcal{H}},$
 where $C^{2}= \max_{0 \leq k \leq \bar{T}} E(|u^{\infty}(k)|^{2}).$ This proves
 the statement for  $\mathfrak{a}_{0}.$
 Since $0 \leq \mathfrak{b}_{0}(\eta)$ $=E(( U(\infty,\eta))^{2})$ $-(E(U(\infty,\eta)))^{2}$
 $\leq  \mathfrak{a}_{0}(\eta),$ it also follows that $\mathfrak{b}_{0}$ is bounded.  \\

\noindent
\textbf{Proof of Theorem \ref{ThA.4}}
We first prove the statement for $\mathfrak{a}.$ Because of continuity of
 $\mathfrak{a},$ which follows from Schwarz inequality and inequalities (\ref{EqA.4}), it is enough to prove inequality
 (\ref{EqA.4.1}) for  $\eta \in \mathcal{D}(\mathfrak{a}_{0})=\mathcal{E}_{\bar{T}}(\mathbb{R}^{N},\mathcal{A}).$
 In this proof we use the notation $X_{k}=\eta (k) \cdot u^{\infty}(k),$ and $Y_{k}=\sum_{0 \leq l \leq k} X_{l},$
 for $0 \leq k \leq \bar{T}.$ By the definition of $\mathfrak{a}_{0}$ (see before (\ref{EqA.4N})) and by the definition
 of $Y_{\bar{T}},$ it follows that $\mathfrak{a}_{0}(\eta)=E(Y_{\bar{T}}^{2}).$
 We shall first prove that
 \begin{equation}
   \mathfrak{a}_{0}(\eta) \geq C E(\sum_{0 \leq k \leq \bar{T}} X_{k}^{2}),
  \label{EqP.4.2.1}
 \end{equation}
 for some $C>0.$
 Let $1 \leq k \leq \bar{T}.$ For $0 \leq l \leq k-1$ we obtain, using that $\eta(l)$
 and $\eta(k)$ are $\mathcal{F}_{k}$ measurable, using the independence of $u^{\infty}(l)$ and $u^{\infty}(k)$
 and using the $\mathcal{F}_{k}$ measurability of $E(X_{k}|\mathcal{F}_{k})$ that
  \begin{equation}
  \begin{split}
 E(X_{l}X_{k})&=E(E((\eta (l) \cdot u^{\infty}(l))(\eta (k) \cdot u^{\infty}(k))|\mathcal{F}_{k})) \\
  =\sum_{i,j}&E(\eta_{i}(l)\eta_{j}(k) E(u^{\infty}_{i}(l)u^{\infty}_{j}(k)|\mathcal{F}_{k})) \\
  =\sum_{i,j}&E(\eta_{i}(l)\eta_{j}(k) E(u^{\infty}_{i}(l)|\mathcal{F}_{k})E(u^{\infty}_{j}(k)|\mathcal{F}_{k})) \\
  =E(&E(\eta (l) \cdot u^{\infty}(l)|\mathcal{F}_{k}) E(\eta (k) \cdot u^{\infty}(k)|\mathcal{F}_{k})) \\
  =E(&E(\eta (l) \cdot u^{\infty}(l) E(\eta (k) \cdot u^{\infty}(k)|\mathcal{F}_{k})|\mathcal{F}_{k})) \\
  =E(&E(\eta (l) \cdot u^{\infty}(l) E(\eta (k) \cdot u^{\infty}(k)|\mathcal{F}_{k})).
  \end{split} \label{EqP.4.2.2}
 \end{equation}
 This proves that if $0 \leq l \leq k-1$ then
 \begin{equation}
   E(X_{l}X_{k})= E(X_{l}E(X_{k}|\mathcal{F}_{k})).
  \label{EqP.4.2.3}
 \end{equation}
 If $1 \leq k \leq \bar{T}$ then $E(Y_{k}^{2})$ $=E((Y_{k-1}+X_{k})^{2})$ $=E(Y_{k-1}^{2}+2Y_{k-1}X_{k}+X_{k}^{2}).$
 Formulas (\ref{EqP.4.2.3}) and $E(X_{k}^{2})$ $=E(E(X_{k}^{2}|\mathcal{F}_{k}))$ now give that
  \begin{equation}
 E(Y_{k}^{2})=E(Y_{k-1}^{2}+2Y_{k-1}E(X_{k}|\mathcal{F}_{k})+E(X_{k}^{2}|\mathcal{F}_{k})),
  \label{EqP.4.2.4}
 \end{equation}
 for $1 \leq k \leq \bar{T}.$ Let $a \in \mathbb{R}^{N},$ $a \neq 0.$ Then
  \begin{equation}
   \begin{split}
 |\sum_{i}& a_{i} E( u^{\infty}_{i}(k)|\mathcal{F}_{k})|=|E(a \cdot u^{\infty}(k))| \\
 &< (E(|a \cdot u^{\infty}(k)|^{2}))^{1/2}=(E(|a \cdot u^{\infty}(k)|^{2}|\mathcal{F}_{k})^{1/2},
  \end{split}   \label{EqP.4.2.5}
 \end{equation}
 where the conditional expectation and the expectation are interchangeable because of hypothesis ($H_{1}$)
 and where there is strict inequality since $a \cdot u^{\infty}(k)$ is a nontrivial random variable according to hypothesis($H_{2}$).
 Let
 \begin{equation}
   \alpha (k) =\max_{|a|=1}|E(a \cdot u^{\infty}(k))|/(E(|a \cdot u^{\infty}(k)|^{2})^{1/2}.
  \label{EqP.4.2.6}
 \end{equation}
 Then $0 \leq \alpha (k) < 1.$ Inequality (\ref{EqP.4.2.5}) and definition (\ref{EqP.4.2.6}) give that
 \begin{equation}
   |E(X_{k}|\mathcal{F}_{k})| \leq \alpha (k) (E(|X_{k}|^{2})^{1/2}.
  \label{EqP.4.2.7}
 \end{equation}
 It follows from formula (\ref{EqP.4.2.4}) and inequality (\ref{EqP.4.2.7}) that 
   \begin{equation*}
   \begin{split}
  E(Y_{k}^{2})& \geq E(Y_{k-1}^{2}-2 |Y_{k-1}|(E(|X_{k}|^{2}|\mathcal{F}_{k}))^{1/2}+E(X_{k}^{2}|\mathcal{F}_{k})) \\
   &= (1- \alpha (k)) E(Y_{k-1}^{2}-2 |Y_{k-1}|(E(|X_{k}|^{2} +E(X_{k}^{2}|\mathcal{F}_{k})) \\
   &+ \alpha (k) E(( |Y_{k-1}|-(E(|X_{k}|^{2}|\mathcal{F}_{k}))^{1/2})^{2}) \\
  &\geq (1- \alpha (k)) E(Y_{k-1}^{2} +E(X_{k}^{2}|\mathcal{F}_{k})).
  \end{split}
 \end{equation*}
 This proves that there exists $C_{k} >0,$ (independent of $\eta $) such that
 \begin{equation}
   E(Y_{k}^{2}) \geq C_{k} E(Y_{k-1}^{2}+X_{k}^{2}),
  \label{EqP.4.2.8}
 \end{equation}
 for $1\leq k \leq \bar{T}.$
 A finite iteration from $k=\bar{T}$ to $k=1$ of this inequality and using that $Y_{0}=X_{0},$
 then proves inequality (\ref{EqP.4.2.1}).

Since $\eta_{j}(k)$ is $\mathcal{F}_{k}$ measurable it follows that
   \begin{equation*}
   \begin{split}
   E(&\sum_{0 \leq k \leq \bar{T}} X_{k}^{2})
     =\sum_{0 \leq k \leq \bar{T}} \sum_{i,j} E(\eta_{i}(k)\eta_{j}(k)  E(u_{i}^{\infty}(k) u_{j}^{\infty}(k)|\mathcal{F}_{k})) \\
     \geq &\sum_{0 \leq k \leq \bar{T} } \sum_{i,j} E(\eta_{i}(k)\eta_{j}(k)  E((u_{i}^{\infty}(k)-E(u_{i}^{\infty}(k))) (u_{j}^{\infty}(k)-E(u_{i}^{\infty}(k)))|\mathcal{F}_{k})), \\
  \end{split}
 \end{equation*}
 where we have used hypothesis ($H_{1}$) to deduce the inequality.
 Then, according to ($H_{2}$):
\begin{equation}
  E(\sum_{0 \leq k \leq \bar{T}} X_{k}^{2}) \geq \sum_{k \geq 0 }
 E(\sum_{i,j} c_{ij}(k) \eta_{i}(k) \eta_{j}(k)).
 \label{EqP.4.2.9}
\end{equation}
This inequality and (\ref{EqP.4.2.1}) prove that  $\mathfrak{a}_{0}(\eta)$
$\geq \sum_{k \geq 0 } c(k)^{2} E( |\eta(k)|^{2})   \geq c^{2} \| \eta \|^{2},$
where $c(k)^{2}$ and $c^{2}$ are some strictly positive numbers. This proves the
inequality for $\mathfrak{a}_{0}(\eta)$ in (\ref{EqA.4.1}).

 To prove the statement for $\mathfrak{b},$ in (\ref{EqA.4.1}), we first show  that if $\eta \in \mathcal{H}$ and $\eta \neq 0$ then  the random variable $ U(\infty, \eta)$ is not constant (a.e.).
 Given $\eta \neq 0,$ let $M$  be the largest integer such that $0 \leq M \leq  \bar{T}$ and $\eta(M) \neq 0.$ If $M=0,$ let $\mathcal{G}$ $=\mathcal{F}_{0}.$ If $M>0,$ let $\mathcal{G}$ be the smallest $\sigma$--algebra, such that $\mathcal{F}_{M}$ $\subset \mathcal{G}$ and such that $u^{\infty}(k)$ is $\mathcal{G}$--measurable for $0 \leq k \leq M-1.$
 By hypotheses ($H_{1}$) and ($H_{3}$), $\mathcal{G}$ and  $u^{\infty}(M)$ are independent, and $X_{k}$ is $\mathcal{G}$--measurable for $0 \leq k \leq M-1.$
 We remind that, if an absolutely integrable random-variable $Z$ is independent of $\mathcal{G},$ then $E(Z | \mathcal{G})= E(Z).$ Using hypothesis ($H_{1}$), ($H_{2}$) and ($H_{3}$) it then follows  for the conditional variance that
  \begin{equation}
   \begin{split}
  E((U(\infty, \eta)&-E( U(\infty, \eta)|\mathcal{G}))^{2}|\mathcal{G})
                  =E((X_{M}-E(X_{M} |\mathcal{G}))^{2}|\mathcal{G}) \\
     =\sum_{i,j}  \eta_{i}(M) & \eta_{j}(M)
        E((u_{i}^{\infty}(M)-E(u_{i}^{\infty}(M)|\mathcal{G}))
           (u_{j}^{\infty}(M)-E(u_{i}^{\infty}(M)|\mathcal{G}))|\mathcal{G}) \\
     =\sum_{i,j}& \eta_{i}(M)\eta_{j}(M)
        E((u_{i}^{\infty}(M)-E(u_{i}^{\infty}(M))) (u_{j}^{\infty}
                                                (M)-E(u_{i}^{\infty}(M)))) \\
      =\sum_{i,j}& \eta_{i}(M)\eta_{j}(M) c_{ij} \geq c^{2} | \eta(M) |^{2}, \\
    \end{split} \label{EqP.4.2.10}
 \end{equation}
 where $c^{2}>0.$
 If $ U(\infty, \eta)$ is a constant, then
$E((U(\infty, \eta)-E( U(\infty, \eta)|\mathcal{G}))^{2}|\mathcal{G})$ $=0,$
which is in contradiction with inequality (\ref{EqP.4.2.10}).
 This proves that $ U(\infty, \eta)$ is not constant, if $\eta \neq 0.$

 Since $ U(\infty, \eta)$ is not constant, if $\eta \neq 0,$ it follows that
 $E((U(\infty, \eta)-E( U(\infty, \eta)))^{2})>0,$ if $\eta \neq 0.$ This proves
 that $\mathfrak{b}(\eta) >0,$ if $\eta \neq 0.$ 

 Let $c=\inf_{\mathfrak{a}(\eta)=1} \mathfrak{b}(\eta).$ Suppose that $c=0.$ Let
 $\{ \eta^{n} \}_{n \geq 1}$
 be a sequence such that $\mathfrak{a}(\eta^{n})=1$ and $\lim_{n \rightarrow \infty } \mathfrak{b}(\eta^{n})=0.$
 Since $\mathfrak{b}(\eta)=\mathfrak{a}(\eta)-(E( U(\infty, \eta)))^{2},$ it follows that
 $(E( U(\infty, \eta^{n})))^{2} \rightarrow 1.$ There is no restriction to suppose
 that $1$ is an accumulation point of $\{ E( U(\infty, \eta^{n}))\},$ since the set $\{ \eta \in \mathcal{H} \, | \, \mathfrak{a}(\eta)=1 \}$ is invariant
 under the transformation $\eta \rightarrow -\eta.$ By selecting a subsequence and
 changing the enumeration, we can suppose that $E( U(\infty, \eta^{n})) \rightarrow 1.$
 Let $\mathcal{H}_{\mathfrak{a}},$ be the Hilbert space defined on $\mathcal{H},$ as
 a linear space, by the scalar product $( \quad , \quad )_{\mathfrak{a}},$ where
 $( \eta , \xi )_{\mathfrak{a}}$ $=( \eta ,A \xi )_{\mathcal{H}}.$  The norms $\mathfrak{a}$
 and $ \| \quad \|^{2}_{\mathcal{H}}$ are equivalent, due to inequalities (\ref{EqA.4})
 and (\ref{EqA.4.1}) for $\mathfrak{a}_{0}$ and $\mathfrak{a},$ so $\mathcal{H}$ and  $\mathcal{H}_{\mathfrak{a}}$
 are identical as TVS's. Since the unit ball of $\mathcal{H}_{\mathfrak{a}}$ is weakly
 compact, ($\mathcal{H}$ is reflexive and Banach-Alaoglu theorem, c.f. \cite{R1}), it  follows
 that $\{ \eta^{n} \}_{n \geq 1}$ has a weakly convergent subsequence. Once more, by selecting
 this  subsequence and  changing the enumeration, we can suppose that
  $\{ \eta^{n} \}_{n \geq 1}$ is weakly convergent
 to an element $\xi \in \mathcal{H}.$ By the definition of weak convergence it
 then follows that $1=\lim_{n \rightarrow} E( U(\infty, \eta^{n}))$ $=E( U(\infty, \xi)).$
 This proves that $\xi \neq 0.$ Moreover, since $\xi$ is an element of the unit
 ball of $\mathcal{H}_{\mathfrak{a}}$ and since $(E( U(\infty, \xi)))^{2} \leq \mathfrak{a}(\xi),$
 it follows that $\mathfrak{a}(\xi)=1.$ Therefore $\mathfrak{b}(\xi)=0,$ which is in contradiction
 with the already proved fact that $\mathfrak{b}(\xi)>0,$ for $\xi \neq 0$ and $\xi \in \mathcal{H}.$
 This proves that $c>0.$ Consequently $\mathfrak{b}(\eta) \geq c \mathfrak{a}(\eta),$
 for $\eta \in \mathcal{H}.$ The statement for $\mathfrak{b}$ in (\ref{EqA.4.1}) now
 follows from the statement for $\mathfrak{a}$ in (\ref{EqA.4.1}). \\

\noindent
\textbf{Proof of Theorem \ref{ThA.5}}
Let $\mathfrak{c}$ $=\mathfrak{b}$ (resp.$\mathfrak{a}$) and
$\mathcal{C}$ $=\mathcal{C}_{1}$ (resp. $\mathcal{C}_{2}$). $\mathcal{C}$ is then
closed and convex in $\mathcal{H}.$ Due to Lemma \ref{LemA.3} and Theorem \ref{ThA.4},
the norm in $\mathcal{H}$ is equivalent to the norm given by $(\mathfrak{c})^{1/2}.$
Let $\mathcal{H}_{\mathfrak{c}}$ be the Hilbertspace with scalar product
$\mathfrak{c}.$ Then $\mathcal{C}$ is closed and convex also in $\mathcal{H}_{\mathfrak{c}}.$
The set $\mathcal{C}$ therefore has a unique point $\hat{\eta},$ minimizing its
distance (in $\mathcal{H}_{\mathfrak{c}}$) to the the origin (c.f. Theorem 12.3 of \cite{R1}).  \\

\noindent
\textbf{Proof of Lemma \ref{LemA.4}}
To prove statement (i), let $\hat{\eta} \in \mathcal{C}_{0}$ be a
 solution of equation (\ref{EqA.0}) and let $E(U(\infty, \hat{\eta}))=e.$
 Then $E(( U(\infty, \hat{\eta}) -E(U(\infty, \hat{\eta})))^{2}) = \sigma^{2},$
 as noted in the paragraph before equation (\ref{EqA.2N}). Let
 $\inf_{\eta \in \mathcal{C}_{1}}    E(( U(\infty, \eta)$ $ -E(U(\infty, \eta)))^{2})
 =\Sigma^{2}.$ Obviously, according to (\ref{EqA.2N}),
 $\Sigma^{2} \leq \sigma^{2},$ since $\hat{\eta} \in \mathcal{C}_{1}.$ Let
 $\Sigma^{2} < \sigma^{2}.$ Then there exists $\eta' \in \mathcal{C}_{1},$ such
 that $ E(( U(\infty, \eta') -E(U(\infty, \eta')))^{2})=(\Sigma')^{2},$ where
 $(\Sigma')^{2} < \sigma^{2}.$
 Then $\eta''= (\sigma / \Sigma')\eta' \in \mathcal{C}_{0}$ and
 $E(U(\infty, \eta''))=(\sigma / \Sigma')e > e,$ which is a contradiction since
 $\hat{\eta}$ is a solution of equation (\ref{EqA.0}). Therefore
 $\sigma = \Sigma$ and $\hat{\eta} \in \mathcal{C}_{1}$ is a solution of
 $E(( U(\infty, \hat{\eta}) -E(U(\infty, \hat{\eta})))^{2})=\Sigma^{2}.$
 To prove statement (ii), let $\hat{\eta} \in \mathcal{C}_{1}$ be a solution of
 equation (\ref{EqA.2N}). Let
 $\sigma^{2}=E(( U(\infty, \hat{\eta}) -E(U(\infty, \hat{\eta})))^{2}),$ let
 $m=E(U(\infty, \hat{\eta}))$ and let
 $M= \sup_{\eta \in \mathcal{C}_{0}} E(U(\infty, \eta)).$
 Obviously $m \leq M$ according to (\ref{EqA.0}), since
 $\hat{\eta} \in \mathcal{C}_{0}.$ Let $m < M.$  Then there exists
 $\eta' \in \mathcal{C}_{0},$ such that
 $\sigma^{2}= E(( U(\infty, \eta') -E(U(\infty, \eta')))^{2})$ and
 $E( U(\infty, \eta'))=M',$ where $m<M'.$ Since $e \leq m < M'$ it follows that
 $\eta' \in \mathcal{C}_{1}.$ This shows that equation (\ref{EqA.2N}) has at least
 two distinct solutions $ \hat{\eta}$ and $ \eta'.$ This is in contradiction,
 with the fact  that the solution of equation
 (\ref{EqA.2N}) is unique, according to Theorem \ref{ThA.5}. Therefore $m=M,$
 which proves that $ \hat{\eta}$ is a solution of (\ref{EqA.0}).    \\

\noindent
\textbf{Proof of Theorem \ref{ThA.7}}
Since the two cases are so similar, we only prove the statements in
 case of $\mathfrak{a}.$ Let $\mathcal{C}_{2}$ be non-empty and let $\hat{\eta}$ be the unique
 solution of equation (\ref{EqA.12}). If $\eta \in \mathcal{C}_{2},$ $t>0$ and  $\eta \neq \hat{\eta},$
 then $\mathfrak{a}(\hat{\eta}+t(\eta-\hat{\eta}))-\mathfrak{a}(\hat{\eta})>0,$ since
 $\hat{\eta}$ is unique and $\mathcal{C}_{2}$ is convex. This gives that
 $(D\mathfrak{a})(\hat{\eta};\epsilon)$ $=(A\hat{\eta},\epsilon) \geq 0,$ for
 $\epsilon \in T_{\mathcal{C}_{2}}(\hat{\eta}),$
 where $T_{\mathcal{C}_{2}}(\hat{\eta})$ is the closed tangent cone of the convex set $\mathcal{C}_{2}$
 at $\hat{\eta}.$ Let $J_{+}(\hat{\eta})$ be the set of all $\epsilon \in \mathcal{H},$ such
 that $\epsilon_{i}(k) \geq 0$ and $(\epsilon_{i}(k))(\omega) = 0$ a.e. $\omega \in supp \, \hat{\eta}_{i}(k),$
 for $ 0 \leq k \leq \bar{T}$ and $1 \leq i \leq N.$ Let $H(\hat{\eta})$ be the
 closed subspace (of $\mathcal{H}$) of elements  $\eta \in \mathcal{H},$ such that
 $supp \, \eta_{i}(k) \subset \, supp \, \hat{\eta}_{i}(k),$ for $ 0 \leq k \leq \bar{T}$ and $1 \leq i \leq N.$
 Then $H(\hat{\eta})$ $=(J_{+}(\hat{\eta}))^{\bot},$ the subset of vectors in $\mathcal{H},$
 orthogonal to every vector in $J_{+}(\hat{\eta}).$ Let $J(\hat{\eta})$ be the set of all
 vectors $\epsilon \in \mathcal{H},$ such that $(\epsilon, \eta)_{\mathcal{H}} \geq 0,$
 for all $\eta \in J_{+}(\hat{\eta}).$ Then $J(\hat{\eta})$ $=J_{+}(\hat{\eta})$
 $+H(\hat{\eta}).$ A vector $\epsilon \in T_{\mathcal{C}_{2}}(\hat{\eta}),$ if and only if
  $\epsilon \in J(\hat{\eta}),$  $(m,\epsilon)_{\mathcal{H}}=0,$ and $(l_{t},\epsilon)_{\mathcal{H}} \geq 0,$
 for $t \in L,$ where
 $L=\{t \in \mathbb{N} \, | \, 0 \leq t \leq \bar{T}+T-1, E((\Delta U)(t+1,\hat{\eta}))- c(t)E(K(t,\hat{\eta} ))=0 \}.$
  Hence
 \begin{equation}
   T_{\mathcal{C}_{2}}(\hat{\eta})=\bigcap_{\eta \in I(\hat{\eta})}
    \{ \epsilon \in \mathcal{H} \, | \, (\eta, \epsilon)_{\mathcal{H}} \geq 0 \},
  \label{EqP.7.1}
 \end{equation}
 where $I(\hat{\eta})$ $=J_{+}(\hat{\eta})$ $\cup \{ m, -m \} \cup \{l_{t} \, | \, t \in L \}.$
 This is equivalent to
 \begin{equation}
   T_{\mathcal{C}_{2}}(\hat{\eta})=\bigcap_{\eta \in V_{+}(\hat{\eta})}
    \{ \epsilon \in \mathcal{H} \, | \, (\eta, \epsilon)_{\mathcal{H}} \geq 0 \},
  \label{EqP.7.2}
 \end{equation}
 where $V_{+}(\hat{\eta})$ is the closed convex cone generated by $I(\hat{\eta}).$
 This shows that $T_{\mathcal{C}_{2}}(\hat{\eta})$ $=(V_{+}(\hat{\eta}))^{\circ},$
 the  (positive) polar of $V_{+}(\hat{\eta}).$ Since $(A\hat{\eta},\epsilon) \geq 0,$
 for all $\epsilon \in T_{\mathcal{C}_{2}}(\hat{\eta}),$ it follows by definition,
 that $A\hat{\eta} \in (T_{\mathcal{C}_{2}}(\hat{\eta}))^{\circ},$ the (positive) polar of
 $T_{\mathcal{C}_{2}}(\hat{\eta}).$ But $((V_{+}(\hat{\eta}))^{\circ})^{\circ}$
 $=V_{+}(\hat{\eta}),$ because $V_{+}(\hat{\eta})$ is a closed convex cone. This
 proves that $A\hat{\eta} \in V_{+}(\hat{\eta}),$ so there are $\mu  \in \mathbb{R},$  $\lambda_{t} \geq 0,$  for $t \in L,$
 and $\nu \in J_{+}(\hat{\eta}),$ such that $A\hat{\eta}$ $=\mu \, m$ $+\sum_{ t \in L} \lambda_{t} l_{t}$ $+\nu.$
 Setting $\lambda_{t}=0,$ for $0 \leq t \leq \bar{T}+T-1$ and $t \notin L,$ we obtain
 formula (\ref{EqA.16}), since $A$ has an inverse defined on $\mathcal{H}.$ Moreover
 conditions (\ref{EqA.13.1}), (\ref{EqA.14}) and (\ref{EqA.15}) are also satisfied. \\

\subsection{Proofs of results in Appendix \ref{Operator C}}
\label{Proof Operator C}
\textbf{Proof of Proposition \ref{PropA.1}}
The equation $(B(k,k)-\lambda) \eta (k) =\xi (k)$ give, according to
formula (\ref{EqA.18}), formula (\ref{EqA.17}), formula (\ref{EqA.19})
and hypothesis ($H_{1}$), that
 $(B(k,k)-\lambda) \eta (k) = (M^{\mathfrak{a}}(k)-\lambda) \eta (k)
    - E(u^{\infty}(k)) (E(u^{\infty}(k)) \cdot E(\eta(k))) = \xi(k),$
for $0 \leq k \leq N.$ Since
$M_{ij}^{\mathfrak{b}}(k) = M_{ij}^{\mathfrak{a}}(k)
    - E(u_{i}^{\infty}(k)) E(u_{j}^{\infty}(k)),$
it follows that
 $(M^{\mathfrak{a}}(k)-\lambda) (\eta (k) -E( \eta (k)))
 +(M^{\mathfrak{b}}(k)-\lambda)E( \eta (k)) = \xi(k).$
Taking the expectation of the two members
of this expression, we first obtain that
 \[(M^{\mathfrak{b}}(k)-\lambda) E( \eta (k)) = E(\xi(k))\]
and then that
 \[(M^{\mathfrak{a}}(k)-\lambda)(\eta (k) -E( \eta (k))) =\xi(k)- E(\xi(k)).\]
The matrices $M^{\mathfrak{a}}(k)-\lambda$ and $M^{\mathfrak{b}}(k)-\lambda$ are invertible
according to the hypothesis that $\lambda \notin \sigma (M^{\mathfrak{a}}(k))$ and
$\lambda \notin \sigma (M^{\mathfrak{b}}(k)).$
The expression (\ref{EqA.20}) of $(B(k,k)- \lambda )^{-1}$ now follows. \\

\noindent
\textbf{Proof of Proposition \ref{PropA.2}}
Using formula (\ref{EqA.23.2}) it follows that
\begin{equation}
   g_{n}=\sum_{n+1 \leq r \leq \bar{T}} d_{r} f_{r}^{2}
 \label{EqP.9.1}
\end{equation}
and using formulas (\ref{EqA.22.3}) and (\ref{EqA.23.3}) it follows that
\begin{equation}
   D_{n}(k)=M^{\mathfrak{a}}(k)
 -\sum_{n+1 \leq r \leq \bar{T}} d_{r} f_{r}^{2} N_{r}^{\mathfrak{a}}(k) - \lambda,
 \label{EqP.9.2}
\end{equation}
where  $0 \leq k \leq n \leq  \bar{T}.$

If $n=\bar{T},$ then formula (\ref{EqA.25.2}) follows from (\ref{EqA.20}). Let
$0 \leq n < \bar{T}$ and let $\xi(k) \in \mathcal{H}_{k}.$ The equation
$B^{n}(k,k) \eta(k)= \xi(k)$ and formula (\ref{EqA.24.7}) give that
\begin{equation}
 D_{n}(k) \eta(k)=\xi(k)+
   (1-g_{n}) E(u^{\infty}(k)) (E(u^{\infty}(k)) \cdot E(\eta(k))),
 \label{EqP.9.3}
\end{equation}
where $\eta(k) \in \mathcal{H}_{k}$ is unknown and $0 \leq  k \leq n.$
Taking the expectation on both sides we obtain
\[
D_{n}(k) E(\eta(k))=E(\xi(k))
  + (1-g_{n}) E(u^{\infty}(k)) (E(u^{\infty}(k)) \cdot E(\eta(k))),\]
where $0 \leq  k \leq n.$ Substitution of (\ref{EqP.9.1}) and (\ref{EqP.9.2})
into this expression gives that
\begin{equation}
  \begin{split}
 (M^{\mathfrak{a}}(k)-\lambda)E( \eta(k))
                   -&E(u^{\infty}(k)) (E(u^{\infty}(k)) \cdot E(\eta(k))) \\
      -\sum_{n+1 \leq r \leq \bar{T}} d_{r} f_{r}^{2} (N_{r}^{\mathfrak{a}}&(k)E(\eta(k))
   -E(u^{\infty}(k)) (E(u^{\infty}(k)) \cdot E(\eta(k)))) \\
    =&E(\xi(k)),
 \end{split}  \label{EqP.9.4}
\end{equation}
for $0 \leq  k \leq n.$ It follows from (\ref{EqA.22.2}) that
\[
   (N_{n}^{\mathfrak{b}}(k))_{ij}=(N_{n}^{\mathfrak{a}}(k))_{ij}
      -E(u^{\infty}_{i}(k))E(u^{\infty}_{j}(k)), \]
which together with (\ref{EqP.9.4}) and the definition of $M^{\mathfrak{b}}(k)$
(see below formula (\ref{EqA.18})) give that
\begin{equation}
 (M^{\mathfrak{b}}(k)
-\sum_{n+1 \leq r \leq \bar{T}}d_{r} f_{r}^{2} N_{r}^{\mathfrak{b}}(k) -\lambda)
 E(\eta(k)) =E(\xi(k)),
 \label{EqP.9.5}
\end{equation}
for $0 \leq  k \leq n.$ According to hypothesis, $\lambda$ is in the resolventset
of
$M^{\mathfrak{b}}(k)
     -\sum_{n+1 \leq r \leq \bar{T}}d_{r} f_{r}^{2} N_{r}^{\mathfrak{b}}(k),$
which proves that
\begin{equation}
   E(\eta(k))=
 (M^{\mathfrak{b}}(k)
-\sum_{n+1 \leq r \leq \bar{T}}d_{r} f_{r}^{2} N_{r}^{\mathfrak{b}}(k) -\lambda)^{-1}
 E(\xi(k)),
 \label{EqP.9.6}
\end{equation}
for $0 \leq  k \leq n.$

Equation (\ref{EqP.9.3}) and the following equation give that
\[ D_{n}(k)(\eta(k) - E(\eta(k)))=\xi(k)-E(\xi(k)).\]
Since $\lambda$ is also in the resolventset
of
$M^{\mathfrak{a}}(k)
     -\sum_{n+1 \leq r \leq \bar{T}}d_{r} f_{r}^{2} N_{r}^{\mathfrak{a}}(k),$
according to hypothesis, it follows by expression (\ref{EqP.9.2})
of $D_{n}(k)$ that
\begin{equation}
   \eta(k) - E(\eta(k))=(M^{\mathfrak{a}}(k)
 -\sum_{n+1 \leq r \leq \bar{T}} d_{r} f_{r}^{2} N_{r}^{\mathfrak{a}}(k) - \lambda)^{-1}
  ( \xi(k)-E(\xi(k))),
  \label{EqP.9.7}
 \end{equation}
for $0 \leq  k \leq n.$ Formulas (\ref{EqP.9.6}) and (\ref{EqP.9.7}) prove expression
(\ref{EqA.25.2}). \\

\noindent
\textbf{Proof of Proposition \ref{PropA.3}}
We shall first find a simplified expression of certain products, of
operators $A(k,l): \mathcal{H}_{l} \rightarrow \mathcal{H}_{k}$ given by
(\ref{EqA.17}). Let $M$ be a real symmetric $N \times N$ matrix and let
\begin{equation}
   I^{n}(k,l)=A(k,n)MA(n,l),
 \label{EqP.10.1}
\end{equation}
for $0 \leq k < n,$ $0 \leq l < n,$ and $0 \leq n \leq \bar{T}.$ The operator
$I^{n}(k,l): \mathcal{H}_{l} \rightarrow \mathcal{H}_{k}$ is linear and
continuous. Definition (\ref{EqA.17}) and formula (\ref{EqA.19}) give that
\begin{equation}
  \begin{split}
   I^{n}(k,l)&\eta(l)=E(u^{\infty}(k)(u^{\infty}(n) \cdot (MA(n,l)\eta(l))) | \mathcal{F}_{k}) \\
      =E&(u^{\infty}(k)(u^{\infty}(n) \cdot E(Mu^{\infty}(n) (u^{\infty}(l) \cdot \eta(l)) | \mathcal{F}_{n})) | \mathcal{F}_{k}),
 \end{split} \label{EqP.10.2}
\end{equation}
for $\eta(l) \in \mathcal{H}_{l}.$
Since $\eta(l)$ is $\mathcal{F}_{n}$-measurable, it follows from hypotheses
($H_{1}$) and ($H_{3}$) that
\begin{equation}
  \begin{split}
   I^{n}(k,l)\eta(l)& \\
      =E(u^{\infty}(k&)(u^{\infty}(n) \cdot E(Mu^{\infty}(n))) (E(u^{\infty}(l)| \mathcal{F}_{n}) \cdot \eta(l)) | \mathcal{F}_{k}) \\
      =\sum_{i} E(E&(u^{\infty}(k)(u^{\infty}(n) \cdot E(Mu^{\infty}(n))) E(u^{\infty}_{i}(l)| \mathcal{F}_{n}) \eta_{i}(l) | \mathcal{F}_{k}) \\
      =\sum_{i}(E(u&^{\infty}(n)) \cdot E(Mu^{\infty}(n))) \\
             &E(E(u^{\infty}(k) E(u^{\infty}_{i}(l)| \mathcal{F}_{n}) | \mathcal{F}_{n}) \eta_{i}(l) | \mathcal{F}_{k}) \\
      =(E(u^{\infty}(&n)) \cdot E(Mu^{\infty}(n)))
           \sum_{i}E(E(u^{\infty}(k) | \mathcal{F}_{n}) E(u^{\infty}_{i}(l)| \mathcal{F}_{n}) \eta_{i}(l) | \mathcal{F}_{k}),
 \end{split} \label{EqP.10.3}
\end{equation}
for $0 \leq k < n,$ $0 \leq l < n,$ and $0 \leq n \leq \bar{T}.$

Let $k \neq l.$ The hypothesis ($H_{3}$) and the first and last member of
equality (\ref{EqP.10.3}) give that
\begin{equation}
  \begin{split}
   I^{n}(k,l)\eta(l)&=(E(u^{\infty}(n)) \cdot E(Mu^{\infty}(n))) \\
      &\sum_{i}E(E(u^{\infty}(k) u^{\infty}_{i}(l)| \mathcal{F}_{n}) \eta_{i}(l) | \mathcal{F}_{k}), \quad k \neq l.
 \end{split} \label{EqP.10.4}
\end{equation}
Since $\eta_{i}(l)$ is $\mathcal{F}_{n}$-measurable and $n>k,$ it follows that
\begin{equation}
   I^{n}(k,l)\eta(l)=(E(u^{\infty}(n)) \cdot E(Mu^{\infty}(n))) A(k,l)\eta(l),
       \quad k \neq l,
 \label{EqP.10.5}
\end{equation}
where $0 \leq k < n$ and $0 \leq l < n.$

Let $k = l.$ The first and last member of equality (\ref{EqP.10.3}), the
$\mathcal{F}_{k}$-mea\-su\-ra\-bi\-lity of $\eta(k)$ and hypothesis ($H_{1}$) give that
\begin{equation}
  \begin{split}
   I^{n}(k,k)\eta(k)&=(E(u^{\infty}(n)) \cdot E(Mu^{\infty}(n))) \\
       &\sum_{i}E(E(u^{\infty}(k) | \mathcal{F}_{n}) E(u^{\infty}_{i}(k)| \mathcal{F}_{n})) \eta_{i}(k),
 \end{split} \label{EqP.10.6}
\end{equation}
where $0 \leq k < n.$ With the notation (\ref{EqA.22.1}), formula (\ref{EqP.10.6})
reads
\begin{equation}
   I^{n}(k,k)\eta(k)=(E(u^{\infty}(n)) \cdot E(Mu^{\infty}(n))) N_{n}^{\mathfrak{a}}(k)\eta(k),
 \label{EqP.10.7}
\end{equation}
where $0 \leq k < n.$

Next we shall consider products as in (\ref{EqP.10.1}), but for operators $B(k,l).$
Let $P$ be a linear continuous operator in $\mathcal{H}_{n},$ such that,
if $\eta(k) \in \mathcal{H}_{n}$ and $E(\eta(n))=0,$ then $P\eta(n)=M\eta(n).$
We introduce a linear and continuous operator
$J^{n}(k,l): \mathcal{H}_{l} \rightarrow \mathcal{H}_{k}$ by
\begin{equation}
   J^{n}(k,l)=B(k,n) P B(n,l),
 \label{EqP.10.8}
\end{equation}
for $0 \leq k < n,$ $0 \leq l < n$ and $0 \leq n \leq \bar{T}.$
We note that
\begin{equation}
   E( B(n,l)\eta(l))=0, \quad 0 \leq l < n.
 \label{EqP.10.9}
\end{equation}
In fact, since $\eta(l)$ is $\mathcal{F}_{l}$-measurable, it follows from
hypotheses ($H_{1}$) and ($H_{3}$), that
\begin{equation}
   E(u^{\infty}(l) \cdot \eta(l))=E(u^{\infty}(l)) \cdot E(\eta(l)),
 \label{EqP.10.10}
\end{equation}
and then that
\begin{equation}
  \begin{split}
   &E(B(n,l) \eta(l))=E(E(u^{\infty}(n)(u^{\infty}(l) \cdot \eta(l) -E(u^{\infty}(l)) \cdot E(\eta(l))) | \mathcal{F}_{n})) \\
     &=\sum_{i}E(E(u^{\infty}(n) u^{\infty}_{i}(l) | \mathcal{F}_{l}) \eta_{i}(l)) -E(u^{\infty}(n))(E(u^{\infty}(l)) \cdot E(\eta(l))) \\
     &=\sum_{i}E(u^{\infty}(n)) E(u^{\infty}_{i}(l)) E( \eta_{i}(l)) -E(u^{\infty}(n))(E(u^{\infty}(l)) \cdot E(\eta(l))) \\
     &=0.
 \end{split} \label{EqP.10.11}
\end{equation}
If $E(\eta(n))=0,$ then $B(k,n)\eta(n)=A(k,n)\eta(n),$ according to (\ref{EqA.17}),
(\ref{EqA.18}), (\ref{EqP.10.10}) and hypothesis ($H_{1}$). This gives together
with formulas (\ref{EqP.10.8}) and (\ref{EqP.10.9}) that
\begin{equation}
   J^{n}(k,l)=A(k,n) M B(n,l),
 \label{EqP.10.12}
\end{equation}
for $0 \leq k < n,$ $0 \leq l < n,$ and $0 \leq n \leq \bar{T}.$ Since definition
(\ref{EqA.17}) of $B(k,l)$ and (\ref{EqP.10.10}) give that
\begin{equation}
   B(r,s)\eta(s)=A(r,s)\eta(s)-E(u^{\infty}(r)) E(u^{\infty}(s) \cdot \eta(s)),
 \label{EqP.10.13}
\end{equation}
for $0 \leq r \leq \bar{T}$ and $0 \leq r \leq \bar{T},$ it follows from
formulas (\ref{EqP.10.1}) and (\ref{EqP.10.12}) that
\begin{equation}
   J^{n}(k,l)\eta(l)=I^{n}(k,l)\eta(l)
            -A(k,n) M E(u^{\infty}(n)) E(u^{\infty}(l) \cdot \eta(l)),
 \label{EqP.10.14}
\end{equation}
for $0 \leq k < n,$ $0 \leq l < n$ and $0 \leq n \leq \bar{T}.$
Noting that (use ($H_{1}$) and ($H_{3}$))
\begin{equation}
  \begin{split}
A(k,n) M E(u^{\infty}(n))&=E(u^{\infty}(k) (u^{\infty}(n) \cdot (ME(u^{\infty}(n)))) | \mathcal{F}_{k}) \\
&=E(u^{\infty}(k)) (E(u^{\infty}(n)) \cdot (ME(u^{\infty}(n)))),
 \end{split} \notag
\end{equation}
for $0 \leq k < n,$ it follows from (\ref{EqP.10.14}) that
\begin{equation}
  \begin{split}
   J^{n}(k,l)\eta(l)&=I^{n}(k,l)\eta(l) \\
            -E(u^{\infty}(k))& (E(u^{\infty}(n)) \cdot (ME(u^{\infty}(n)))) (E(u^{\infty}(l)) \cdot E(\eta(l))),
 \end{split}  \label{EqP.10.15}
\end{equation}
for $0 \leq k < n,$ $0 \leq l < n,$  $0 \leq n \leq \bar{T}$ and $\eta(l) \in \mathcal{H}_{l}.$
Formulas (\ref{EqP.10.4}) and (\ref{EqP.10.15}) give that
\begin{equation}
  \begin{split}
J^{n}(k,l)\eta(l)&=(E(u^{\infty}(n)) \cdot (ME(u^{\infty}(n)))) \\
                        &(A(k,l)\eta(l) - E(u^{\infty}(k))(E(u^{\infty}(l)) \cdot E(\eta(l)))),
   \quad k \neq l.
 \end{split} \notag
\end{equation}
Definition (\ref{EqA.17}) of $B(k,l)$ and and formula (\ref{EqP.10.13}) then
give that
\begin{equation}
J^{n}(k,l)\eta(l)=(E(u^{\infty}(n)) \cdot (ME(u^{\infty}(n)))) B(k,l)\eta(l),
   \quad k \neq l,
  \label{EqP.10.16}
\end{equation}
where  $0 \leq k < n,$ $0 \leq l < n,$  $0 \leq n \leq \bar{T}$ and
$\eta(l) \in \mathcal{H}_{l}.$ Formulas (\ref{EqP.10.7}) and (\ref{EqP.10.15})
give that
\begin{equation}
  \begin{split}
J^{n}(k,k)\eta(k)&=(E(u^{\infty}(n)) \cdot (ME(u^{\infty}(n)))) \\
        &(N_{n}^{\mathfrak{a}}(k)\eta(k) - E(u^{\infty}(k))  (E(u^{\infty}(k)) \cdot E(\eta(k)))),
  \end{split}  \label{EqP.10.17}
\end{equation}
where  $0 \leq k < n,$ $0 \leq n \leq \bar{T}$ and
$\eta(k) \in \mathcal{H}_{k}.$

Having established (\ref{EqP.10.5}) (resp. (\ref{EqP.10.16})) and (\ref{EqP.10.7})
(resp. (\ref{EqP.10.17})) in the case of $C=A$ (resp. $B$), we now turn to the
main part of the proof.

For $0 \leq n \leq \bar{T},$ let $P(n)$ be the statement that the following points
(i), (ii) and (iii) are true:
\begin{itemize}
\item i)  $f_{\bar{T}}, \ldots ,f_{n},$ $g_{\bar{T}}, \ldots ,g_{n}$ and the
matrix elements of $D_{s}(k),$ where $0 \leq k \leq s$ and $n \leq s \leq \bar{T},$
are real valued rational functions of $\lambda$ in $\mathbb{R},$ without singularity
in $\mathbb{R} - \sigma_{n+1}^{\mathfrak{a}}.$

If $n < \bar{T},$ then $D_{\bar{T}}(\bar{T}), \ldots ,D_{n+1}(n+1)$ are invertible
matrices, for $\lambda \in \mathbb{R} - \sigma_{n+1}^{\mathfrak{a}},$ and
$d_{0}, \ldots ,d_{\bar{n+1}}$ are real valued rational functions of
$\lambda$ in $\mathbb{R},$ without singularity in
$\mathbb{R} - \sigma_{n+1}^{\mathfrak{a}}.$

(We remind that these functions are defined by (\ref{EqA.22.1}) and
 (\ref{EqA.22.3})-(\ref{EqA.23.3}) and that $\sigma_{n+1}^{\mathfrak{a}}$
is defined by (\ref{EqA.27.1}))
\item ii) If $\lambda \in \mathbb{R} - \sigma_{n+1}^{\mathfrak{c}},$ then
$C^{s}(k,l),$ where $n \leq s \leq \bar{T},$ $0 \leq k \leq \bar{T}$ and
$0 \leq l \leq \bar{T},$ given by (\ref{EqA.24.1})-(\ref{EqA.24.7}),
are continuous linear operators from $\mathcal{H}_{l}$ to $\mathcal{H}_{k}.$
If moreover $n<\bar{T},$ then the operator $C^{s}(s,s),$ where $n+1 \leq s \leq \bar{T},$
have a continuous inverse in $\mathcal{H}_{s}.$
\item iii) To formulate this part of $P(n),$ we introduce linear continuous
operators $S_{\bar{T}}, \ldots ,S_{n}$ in $\mathcal{H},$ defined by $S_{\bar{T}}\eta=\eta$
and, if $n< \bar{T}$ and $n \leq s < \bar{T},$ by
\begin{equation}
  (S_{s}\eta(k))=\eta(k), \quad s+1 \leq  k  \leq \bar{T}
 \label{EqP.10.18}
\end{equation}
and
\begin{equation}
  (S_{s}\eta(k))=\eta(k) -C^{s+1}(k,s+1)(C^{s+1}(s+1,s+1))^{-1}\eta(s+1),
 \label{EqP.10.19}
\end{equation}
for $0 \leq k < s+1.$ (It follows from statement (ii) of $P(n),$ that the operator
$S_{s}: \mathcal{H} \rightarrow \mathcal{H}$ is well defined.)

Let $\lambda \in \mathbb{R} - \sigma_{n+1}^{\mathfrak{c}}.$ Then
$S_{s}: \mathcal{H} \rightarrow \mathcal{H}$ is continuous with continuous
inverse,  for $n \leq s \leq \bar{T},$ and
\begin{equation}
  C^{s}=U_{s}C
 \label{EqP.10.20}
\end{equation}
\begin{equation}
  \xi^{s}=U_{s}\xi,
 \label{EqP.10.21}
\end{equation}
for $n \leq s \leq \bar{T},$ where $U_{s}=S_{s}S_{s+1} \cdots S_{\bar{T}}.$
\end{itemize}

We prove, by finite induction, that $P(\bar{T}), \ldots ,P(0)$ are true. The
statement $P(\bar{T})$ is trivially true. In fact
$\mathbb{R} - \sigma_{\bar{T}+1}^{\mathfrak{a}}=\mathbb{R} - \sigma_{\bar{T}+1}^{\mathfrak{a}}=\mathbb{R}$
(see above (\ref{EqA.27.1})) and $f_{\bar{T}},$ $g_{\bar{T}},$  $D_{\bar{T}}(k),$ where
$0 \leq k \leq \bar{T},$ given by (\ref{EqA.22.3}), are well defined for $\lambda \in \mathbb{R}.$
Statement (i) is then true. Since $C: \mathcal{H} \rightarrow \mathcal{H}$ is
continuous, it follows that $C^{\bar{T}}(k,l): \mathcal{H}_{l} \rightarrow \mathcal{H}_{k},$
$0 \leq k \leq \bar{T},$ $0 \leq l \leq \bar{T},$ given by (\ref{EqA.24.1}) and
(\ref{EqA.24.2}) are continuous. Thus, statement (ii) is true. Statement (iii) is
also true, since $U_{\bar{T}}=S_{\bar{T}}=I,$ where $I$ is the identity in
$\mathcal{H}.$

Let $0< n \leq \bar{T}.$ We shall prove that, if $P(\bar{T}), \ldots ,P(n)$ are
true, then $P(n-1)$ is true. We remind that
$\sigma_{n+1}^{\mathfrak{c}} \subset \sigma_{n}^{\mathfrak{c}},$ (see below (\ref{EqA.27.2})).

We first prove that (i) of $P(n-1)$ is true. If $n=\bar{T},$ then $D_{\bar{T}}(\bar{T}),$
given by (\ref{EqA.22.3}), is invertible for
$\lambda \in \mathbb{R} - \sigma_{\bar{T}}^{\mathfrak{a}}.$ If $n< \bar{T},$ then
$D_{\bar{T}}(\bar{T}), \ldots ,D_{n+1}(n+1)$ are, according to (i) of $P(n),$ invertible
matrices, for $\lambda \in \mathbb{R} - \sigma_{n}^{\mathfrak{a}} \subset  \mathbb{R} - \sigma_{n+1}^{\mathfrak{a}}$
and $D_{n}(n)$ is, according to expression (\ref{EqP.9.2}) of $D_{n}(n)$ and
definition (\ref{EqA.27.1}) of $\sigma_{n}^{\mathfrak{a}},$ also invertible for
$\lambda \in \mathbb{R} - \sigma_{n}^{\mathfrak{a}}.$ Thus,
$D_{\bar{T}}(\bar{T}), \ldots ,D_{n}(n)$ are invertible matrices, for
$\lambda \in \mathbb{R} - \sigma_{n}^{\mathfrak{a}}.$ The matrix elements of the matrices
$D_{\bar{T}}(\bar{T}), \ldots ,D_{n}(n)$ are, according to (i) of $P(n),$
rational functions (of $\lambda$), without singularities in
$\mathbb{R} - \sigma_{n}^{\mathfrak{a}} \, (\subset  \mathbb{R} - \sigma_{n+1}^{\mathfrak{a}}).$
Formulas (\ref{EqA.23}) and (\ref{EqA.23.1}) then show that
$d_{\bar{T}}, \ldots ,d_{n}$ are rational functions, without singularities in
$\mathbb{R} - \sigma_{n}^{\mathfrak{a}}.$ This proves that the second part of
statement (i) of $P(n-1)$ is true.

According to statement (i) of $P(n),$ it follows that
$f_{\bar{T}}, \ldots ,f_{n},$ $g_{\bar{T}}, \ldots ,g_{n}$ and the
matrix elements of $D_{s}(k),$ where $0 \leq k \leq s$ and $n \leq s \leq \bar{T},$
are real valued rational functions of $\lambda$ in $\mathbb{R},$ without singularity
in $\mathbb{R} - \sigma_{n}^{\mathfrak{a}} \subset \mathbb{R} - \sigma_{n+1}^{\mathfrak{a}}.$
Since we already have proved, that $d_{n}$ is a rational function without singularity
in $\mathbb{R} - \sigma_{n}^{\mathfrak{a}},$ it follows from (\ref{EqA.23.2}) and
(\ref{EqA.23.3}), that $f_{n-1},$ $g_{n-1}$ and the matrix elements of
$D_{n-1}(k),$ where $0 \leq k \leq n-1,$
are rational functions of $\lambda$ in $\mathbb{R},$ without singularity
in $\mathbb{R} - \sigma_{n}^{\mathfrak{a}}.$ This completes the proof of statement (i)
of $P(n-1).$

Next we prove statement (ii) of $P(n-1).$ According to (i) of $P(n-1),$ which we
already have proved, $D_{n}(n),$ is invertible, $f_{n-1},$ $g_{n-1},$ are real
numbers and $D_{n-1}(k),$ where $0 \leq k \leq n-1,$ are real matrices, for
$\lambda \in \mathbb{R} - \sigma_{n}^{\mathfrak{c}} \, ( \subset \mathbb{R} - \sigma_{n+1}^{\mathfrak{a}}).$
Let $\lambda \in \mathbb{R} - \sigma_{n}^{\mathfrak{c}}.$
According to (ii) of $P(n),$ $C^{s}(k,l),$ where $n \leq s \leq \bar{T},$ $0 \leq k \leq \bar{T}$
and $0 \leq l \leq \bar{T},$
are continuous linear operators from $\mathcal{H}_{l}$ to $\mathcal{H}_{k}.$
Now it follows from (\ref{EqA.24.3})-(\ref{EqA.24.7}), that
$C^{n-1}(k,l): \mathcal{H}_{l} \rightarrow \mathcal{H}_{k},$ where
$0 \leq k \leq \bar{T}$ and $0 \leq l \leq \bar{T},$ are continuous linear
operators. This proves the first part of statement (ii) of $P(n-1).$

Let $\lambda \in \mathbb{R} - \sigma_{n}^{\mathfrak{c}},$ where $\mathfrak{c}=\mathfrak{a}$
(resp. $\mathfrak{b}$), when $C=A$ (resp. $B$). If $n=\bar{T},$ then
$\sigma_{\bar{T}}^{\mathfrak{a}} = \sigma(M^{\mathfrak{a}}(\bar{T}))$
(resp. $\sigma(M^{\mathfrak{a}}(\bar{T})) \cup \sigma(M^{\mathfrak{b}}(\bar{T}))$)
according to (\ref{EqA.27.1}) (resp. (\ref{EqA.27.2})).
We remind that $D_{\bar{T}}(\bar{T}) = M^{\mathfrak{a}}(\bar{T})-\lambda,$
(see (\ref{EqA.22.3})). According to the inversion formula (\ref{EqA.25.1})
(resp. (\ref{EqA.25.2})) of $A^{\bar{T}}(\bar{T},\bar{T})$
(resp.$B^{\bar{T}}(\bar{T},\bar{T})$) and
its hypothesis, it follows that $A^{\bar{T}}(\bar{T},\bar{T})$
(resp.$B^{\bar{T}}(\bar{T},\bar{T})$) has a continuous inverse in $\mathcal{H}_{\bar{T}}.$
If $n<\bar{T},$ then (ii) of $P(n)$ gives that $C^{s}(s,s)$ has a continuous
inverse in $\mathcal{H}_{s},$ for $n+1 \leq s \leq \bar{T}.$ Since
$\lambda \in \mathbb{R} - \sigma_{n}^{\mathfrak{c}},$ it follows, in the case of
$C=A$ from expression (\ref{EqP.9.2}) of $D_{n}(n)$ and formula (\ref{EqA.25.1})
and in the case of $C=B$ from formula (\ref{EqA.25.2}), that $C^{n}(n,n)$ has
a continuous inverse in $\mathcal{H}_{n}.$ This completes the proof of statement
(ii) of $P(n-1).$

Finally we prove statement (iii) of $P(n-1).$
Let $\lambda \in \mathbb{R} - \sigma_{n}^{\mathfrak{c}}.$
Then $\lambda \in \mathbb{R} - \sigma_{n+1}^{\mathfrak{c}},$ so it follows from
statement (iii) of $P(n),$ that $S_{s}$ is continuous with continuous inverse in
$\mathcal{H},$ for $n \leq s \leq \bar{T}.$ The operators
$C^{n}(k,l): \mathcal{H}_{l} \rightarrow \mathcal{H}_{k},$
$0 \leq k \leq \bar{T}$ and $0 \leq l \leq \bar{T}$ and the inverse of $C^{n}(n,n)$
are continuous according to (the already proved) statement (ii) of $P(n-1).$
Therefore, $S_{n-1}$ given by (\ref{EqP.10.18}) and (\ref{EqP.10.19}) is a continuous
operator in $\mathcal{H}.$ It is invertible, because its inverse is obtained by
formulas (\ref{EqP.10.18}) and (\ref{EqP.10.19}) (in the case of $P(n-1)$)
with the sign changed, of the second term of the right hand side of (\ref{EqP.10.19}).
This proves, that the operators $S_{s}:\mathcal{H} \rightarrow \mathcal{H},$
$n-1 \leq s \leq \bar{T},$ are continuous with continuous inverse.

The operator $C^{s}:\mathcal{H} \rightarrow \mathcal{H},$ satisfy (\ref{EqP.10.20}),
for $n \leq s \leq \bar{T},$ according to statement (iii) of $P(n).$ To establish
(\ref{EqP.10.20}) for $s=n-1,$ it is therefore enough to prove that
$C^{n-1}=S_{n-1}C^{n}.$ Let $C'=S_{n-1}C^{n}.$ It follows from (\ref{EqP.10.18})
and (\ref{EqP.10.19}) that
\begin{equation}
  C'(k,l)=C^{n}(k,l), \quad n \leq  k  \leq \bar{T}, \quad 0 \leq  l  \leq \bar{T},
 \label{EqP.10.22}
\end{equation}
and
\begin{equation}
  C'(k,l)=C^{n}(k,l) -C^{n}(k,n)(C^{n}(n,n))^{-1}C^{n}(n,l),
 \label{EqP.10.23}
\end{equation}
for $0 \leq k < n$ and $0 \leq  l  \leq \bar{T}.$
To prove that $C'(k,l)=C^{n-1}(k,l),$ for $0 \leq  k  \leq \bar{T}$ and
$0 \leq  l  \leq \bar{T},$ we shall consider $(k,l)$ chosen according to a
certain disjoint partition of the set
$\{(k,l) \in \mathbb{N} \times \mathbb{N} \, |\, 0 \leq  k  \leq \bar{T}, \, 
0 \leq  l  \leq \bar{T}\}.$

Formulas (\ref{EqA.24.3}) and (\ref{EqP.10.22}) show that $C'(k,l)=C^{n-1}(k,l),$
for $n \leq  k  \leq \bar{T}$ and $0 \leq  l  \leq \bar{T}.$

Let $0 \leq k < n$
and $0 \leq  l  \leq \bar{T}.$ Formula (\ref{EqP.10.23}) gives that $C'(k,n)=0.$
If $l>n,$ then statement (iii) of $P(n)$ gives that formula (\ref{EqA.24.4}), with
$n+1$ instead of $n,$ is satisfied, so $C^{n}(r,s)=0$ for $0 \leq r < n+1,$
$n+1 \leq s \leq \bar{T}.$ In particular, for the above chosen $k$ and $l,$
$C^{n}(k,l)=0$ and $C^{n}(n,l)=0,$ so $C'(k,l)=0$ according to (\ref{EqP.10.23}).
Formula (\ref{EqA.24.4}) now
gives that $C'(k,l)=C^{n-1}(k,l),$ for $0 \leq k < n$ and $0 \leq  l  \leq \bar{T}.$

Since statement (iii) of $P(n)$ is true, it follows from the expression of $C^{n}(r,s)$
obtained, by expression (\ref{EqA.24.1}) in the case when $n=\bar{T}$ and by expression
(\ref{EqA.24.5}), with $n$ replaced by $n+1,$ in the case when $n<\bar{T}$, that
\begin{equation}
 C^{n}(r,s)= f_{n}C(r,s),
 \label{EqP.10.24}
\end{equation}
where $r \neq s,$ $0 \leq r \leq n$ and $0 \leq s \leq n.$ Formulas (\ref{EqP.10.23})
and (\ref{EqP.10.24}) give that
\begin{equation}
  C'(k,l)=C^{n}(k,l) - f_{n}^{2}C(k,n)(C^{n}(n,n))^{-1}C(n,l),
 \label{EqP.10.25}
\end{equation}
for $0 \leq k < n$ and $0 \leq  l  < n.$
Since $\lambda \in \mathbb{R} - \sigma_{n}^{\mathfrak{c}},$ it follows , in the case of
$C=A$ (resp. $B$), from formula (\ref{EqA.25.1}) (resp. (\ref{EqA.25.2})) and
its hypothesis that $(C^{n}(n,n))^{-1}$ exists and that
$(C^{n}(n,n))^{-1}\eta(n) =(D_{n}(n))^{-1}\eta(n),$ if $\eta(n) \in \mathcal{H}_{n}$
(resp. $\eta(n) \in \mathcal{H}_{n}$ and $E(\eta(n))=0$). The product
$C(k,n)(C^{n}(n,n))^{-1}C(n,l)$ in (\ref{EqP.10.25}) can therefore be simplified,
in the case of $C=A$ (resp. $B$), by using the already derived expressions
(\ref{EqP.10.5}) (resp. (\ref{EqP.10.16})) and (\ref{EqP.10.7}) (resp. (\ref{EqP.10.17}))
of $I^{n}(k,l)$ (resp. $J^{n}(k,l)$), defined by (\ref{EqP.10.1}) (resp. (\ref{EqP.10.8})),
with $M=(D_{n}(n))^{-1}$ (resp. $P=(B^{n}(n,n))^{-1}$ and $M=(D_{n}(n))^{-1}$).

Let $k \neq l,$ $0 \leq k < n$ and $0 \leq  l  < n.$ Formulas (\ref{EqP.10.25}),
(\ref{EqP.10.5}) (resp. (\ref{EqP.10.16})) and expressions (\ref{EqA.23}) and
(\ref{EqA.23.1}) of $d_{n},$ give that
\[
C'(k,l)=C^{n}(k,l) - f_{n}^{2}d_{n}C(k,l). \]
Formula (\ref{EqP.10.24}) gives that
\[
C'(k,l)=f_{n}(1 - f_{n}d_{n})C(k,l). \]
Expression (\ref{EqA.23.2}) of $f_{n-1}$ and expression (\ref{EqA.24.5}) now show
that $C'(k,l)=C^{n-1}(k,l),$ for $k \neq l,$ $0 \leq k < n$ and $0 \leq  l  < n.$

Let $k = l$ and $0 \leq k < n.$ Since $P(n)$ is true, it follows that $C^{n}(k,k)$
is given by (\ref{EqA.24.6}) (resp. (\ref{EqA.24.7})), with $n+1$ instead of $n,$
for $C=A$ (resp. $B$). We first consider the case of $C=A.$ Formulas (\ref{EqP.10.25}),
(\ref{EqA.24.6}) and (\ref{EqP.10.7}) and formulas (\ref{EqA.23}) and
(\ref{EqA.23.1}) of $d_{n},$ give that
\[
  A'(k,k)=D_{n}(k) - f_{n}^{2} d_{n}N_{n}^{\mathfrak{a}}(k)). \]
Definitions (\ref{EqA.23.3}) of $D_{n-1}(k)$ and (\ref{EqA.24.6}) of $A^{n-1}(k,k)$
now show that $A'(k,k)$ $=A^{n-1}(k,k).$ Next we consider the case of $C=B.$ Formulas
(\ref{EqP.10.25}), (\ref{EqA.24.7}) and (\ref{EqP.10.17}) and formulas
(\ref{EqA.23}) and (\ref{EqA.23.1}) of $d_{n},$ give that
\begin{equation}
  \begin{split}
  B'(k,k)\eta(k)&=(D_{n}(k) - f_{n}^{2} d_{n}N_{n}^{\mathfrak{a}}(k))\eta(k) \\
       &-(1-g_{n}-f_{n}^{2} d_{n})
           E(u^{\infty}(k))  (E(u^{\infty}(k)) \cdot E(\eta(k))),
  \end{split}  \notag
\end{equation}
for $\eta(k) \in \mathcal{H}_{k}.$ Definitions (\ref{EqA.23.2}) of $g_{n-1},$
(\ref{EqA.23.3}) of $D_{n-1}(k)$ and (\ref{EqA.24.7}) of $B^{n-1}(k,k)$ now show
that $B'(k,k)=B^{n-1}(k,k).$ This completes the proof of formula (\ref{EqP.10.20}),
for $n-1 \leq s \leq \bar{T}.$

The elements $\xi^{s} \in \mathcal{H}$ satisfy (\ref{EqP.10.20}),
for $n \leq s \leq \bar{T},$ according to statement (iii) of $P(n).$ To establish
(\ref{EqP.10.20}) for $s=n-1,$ it is therefore enough to prove that
$\xi^{n-1}=S_{n-1}\xi^{n}.$
Let $\xi'=S_{n-1}\xi^{n}.$ We remind that $(C^{n}(n,n))^{-1}$ exists (see below
formula (\ref{EqP.10.25})). The definition of $S_{n-1},$ by formulas (\ref{EqP.10.18})
and (\ref{EqP.10.19}), (with $n-1$ instead of $n$) gives that
\[
 \xi'(k)=\xi^{n}(k), \quad n \leq  k  \leq \bar{T} \]
and
\[
\xi'(k)=\xi^{n}(k) -C^{n}(k,n)(C^{n}(n,n))^{-1}\xi^{n}(n),
  \quad   0 \leq k < n. \]
Formulas (\ref{EqA.26}) and (\ref{EqA.27}) defining $\xi^{n-1},$ now show that
$\xi'=\xi^{n-1}.$ This completes the proof of statement (iii) of $P(n-1),$ so
$P(n-1)$ is true.

The equation $C^{\bar{T}} \eta = \xi,$ $\xi \in \mathcal{H}$ (i.e. equation (\ref{EqA.22}))
is, according to $P(n),$ (formulas (\ref{EqP.10.20}) and (\ref{EqP.10.21}) of
statement (iii) of $P(n)$) equivalent to $C^{n} \eta = \xi^{n},$ $0 \leq  n  \leq \bar{T},$
i.e. equation (\ref{EqA.28}). This completes the proof. \\

\noindent
\textbf{Proof of Proposition \ref{PropA.4}}
Let $\lambda \in \mathbb{R} - \sigma_{0}^{\mathfrak{c}}$ and let
$0 \leq n \leq \bar{T}.$ According to the definition of $\sigma_{n}^{\mathfrak{c}}$
(see (\ref{EqA.27.1}) and (\ref{EqA.27.2})), it follows that
$\lambda \in \mathbb{R} - \sigma_{n}^{\mathfrak{c}}.$ Expression (\ref{EqP.9.2}) of
$D_{n}(n),$ formula (\ref{EqA.25.1}) and the definition of $\sigma_{n}^{\mathfrak{a}},$
show that $A^{n}(n,n)$ has a bounded inverse. Formula (\ref{EqA.25.2}) and the
definition of $\sigma_{n}^{\mathfrak{b}},$ show that $B^{n}(n,n)$ has a bounded
inverse. Successive use of (\ref{EqA.24.3}), give that $C^{0}(k,k)=C^{k}(k,k),$
for $0 \leq k \leq \bar{T}.$ Expression (\ref{EqA.30}) now follows from (\ref{EqA.29}).
This completes the proof. \\

\bibliographystyle{amsplain}

\end{document}